\documentclass[10pt]{article}

\usepackage{a4wide}
\usepackage{amssymb}
\usepackage{amsfonts}
\usepackage{amsmath}
\usepackage[all]{xy}

\date{}

\newtheorem{proposition}{Proposition}[section]
\newtheorem{theorem}[proposition]{Theorem}
\newtheorem{lemma}[proposition]{Lemma}

\newtheorem{corollary}[proposition]{Corollary}

\def\der{\partial }

\def\nFM0{{\nu }_{F,M_0}}
\def\nFN0{{\nu }_{F,N_0}}
\def\nGN0{{\nu }_{G,N_0}}

\def\N0{ {\bf N}_0 }

\def\g{\gamma}

\def\ra{\rightarrow}

\def\Xpm{X^{\pm }}

\def\s{\sigma}

\def\l1{{\lambda}_1}

\def\a{\alpha}
\def\a0{ {\alpha }_0}
\def\a1{ {\alpha }_1}

\def\l{\lambda}

\def\nFGM0{{\nu }_{F,G,M_0}}

\def\nFN0{{\nu}_{F,N_0}}

\def\sm{{\sigma}^m}

\def\sm1{{\sigma}^{-1}}

\def\smtp1{{\sigma}^{-t+1}}

\def\S1{S^{-1}}

\def\Xpm1{X^{\pm 1}_1}

\def\sPM1{{\sigma }^{\pm 1}}
\def\sMP1{{\sigma }^{\mp 1 }}

\def\di{{\rm d.ind}}

\def\L{\Lambda}

\def\Ytm1{Y^{t-1}}
\def\Yim1{Y^{i-1}}

\def\CL{{\cal L}}
\def\CM{{\cal M}}
\def\CN{{\cal N}}

\def\ass{{\rm ass}}

\def\supp{{\rm supp }}

\def\ker{ {\rm ker } }

\def\SL2Z{ {\rm SL}_2({\bf Z}) }

\def\th{ \theta }

\def\CL{{\cal L}}

\def\Gp1{ G^{1 , 1 } }
\def\P11{ P^{-1 , 1 } }
\def\Pp1{ P^{1 , 1 } }

\def\th{\theta}

\def\nCLsr{{}^\nu\kern-2pt {\cal L}^{\sigma , \rho  }}
\def\nP{{}^\nu \kern-2pt P}
\def\nL{{}^\nu\kern-2pt L}
\def\nLL{{}^\nu\kern-2pt \Lambda}
\def\nPsr{{}^\nu\kern-2pt P^{\sigma , \rho  }}
\def\nLsr{{}^\nu\kern-2pt L^{\sigma , \rho  }}
\def\nuCL{{}^\nu\kern-2pt  {\cal L}}
\def\nCLsr{{}^\nu\kern-2pt {\cal L}^{\sigma , \rho  }}
\def\nCL1m{{}^\nu\kern-2pt {\cal L}^{-1 , 1  }}
\def\x1nu{x^\frac{1}{\nu}}
\def\xm1nu{x^{-\frac{1}{\nu}}}

\def\CN{{\cal N}}
\def\ra{\rightarrow }

\def\CB{{\cal B}}

\def\CT{{\cal T}}

\def\CC{ {\cal C}}

\def\nAM0{{\nu }_{{\cal A},M_0}}
\def\nAN0{{\nu }_{{\cal A},N_0}}

\def\bR{\overline{R}}

\def\ga{\mathfrak{a}}
\def\gb{\mathfrak{b}}
\def\gc{\mathfrak{c}}

\def\gn{\mathfrak{n}}

\def\gp{\mathfrak{p}}

\def\SL{{\rm SL}}

\def\di!{\frac{\der^i}{i!}}
\def\dik!{\frac{\der^k_i}{k!}}

\def\id{{\rm id}}

\def\Max{{\rm Max}}

\def\N{\mathbb{N}}

\def\0{\overline{0}}
\def\1{\overline{1}}

\def\Ln1{\L_{n,\overline{1}}}

\def\a1{a_{\overline{1}}}

\def\bs{\overline{s}}

\def\S{\Sigma}

\def\vn1{\overrightarrow{n-1}}

\def\im{{\rm im}}

\def\Min{{\rm Min}}

\def\mJ{\mathbb{J}}
\def\mI{\mathbb{I}}

\def\ann{{\rm ann}}
\def\lann{{\rm l.ann}}

\def\K1{{\rm K}_1}

\def\hmI1{\widehat{\mI_1}}
\def\tmI1{\widetilde{\mI_1}}
\def\tmJ1{\widetilde{\mJ_1}}
\def\hB1{\widehat{B_1}}
\def\hCB1{\widehat{\CB_1}}

\def\bS{\overline{S}}

\def\Den{{\rm Den}}

\def\Ore{{\rm Ore}}

\def\Den{{\rm Den}}

\def\Loc{{\rm Loc}}

\def\Ass{{\rm Ass}}

\def\maxDen{{\rm max.Den}}
\def\maxAss{{\rm max.Ass}}
\def\maxLoc{{\rm max.Loc}}
\def\llrad{{\rm l.lrad}}

\def\assmaxDen{{\rm ass.max.Den}}

\def\br{\overline{r}}

\def\bs{\overline{s}}

\def\ga{\mathfrak{a}}

\def\gll{\mathfrak{l}}

\def\pil{\pi_{\gll}}
\def\CCRp{\CC_{R/\gp}}
\def\pilS{\pi_{\gll, S}}

\begin{document}

\author{V. V. \  Bavula 
}

\title{New criteria for a ring to have a semisimple left quotient ring}

\maketitle

\begin{abstract}
Goldie's Theorem (1960), which is one of the most important results in Ring Theory,  is a criterion for a ring to have a semisimple left quotient ring. The aim of the paper is to give four  new criteria (using a completely different approach and new ideas). The first one is based on the recent fact that for an {\em arbitrary} ring $R$ the set $\CM$ of {\em maximal} left denominator sets of $R$ is a non-empty set \cite{larglquot}:

$\noindent $

{\bf Theorem (The First Criterion)}. {\em A ring $R$ has a semisimple left quotient ring $Q$ iff $\CM $ is a finite set, $\bigcap_{S\in \CM } \ass (S) =0$ and,  for each $S\in \CM$, the ring  $S^{-1}R$ is a simple left Artinian ring. In this case, $Q\simeq \prod_{S\in \CM} S^{-1}R$.  }

$\noindent $

 The Second Criterion is given via the minimal primes of $R$ and goes further then the First one 
  in the sense that it describes  explicitly the maximal left denominator sets $S$ via the minimal primes of $R$. The Third Criterion is close to Goldie's Criterion but it is easier to check in applications (basically, it reduces Goldie's Theorem to the prime case). The Fourth Criterion is given via certain left denominator sets.

$\noindent $



$${\bf Contents}$$
\begin{enumerate}
\item Introduction.
\item Preliminaries.
\item  The First Criterion (via the maximal left denominator sets).
 \item The Second  Criterion  (via the minimal primes).
 \item The Third Criterion (in the spirit of Goldie-Lesieur-Croisot).
 \item The Fourth Criterion (via certain left denominator sets).
 \item Left denominator sets of finite direct products of rings.
     \item Criterion for $R/\gll_R$ to have a semisimple left quotient  ring.
\end{enumerate}
\end{abstract}


\section{Introduction}

In this paper, module means a left module, and the following notation is fixed:
\begin{itemize}
\item  $R$ is a ring with 1, $\gn =\gn_R$ is its prime radical and $\Min (R)$ is the set of minimal primes of $R$;
\item   $\CC = \CC_R$  is the set of regular elements of the ring $R$ (i.e. $\CC$ is the set of non-zero-divisors of the ring $R$);
\item   $Q=Q_{l,cl}(R):= \CC^{-1}R$ is the {\em left quotient ring}  (the {\em classical left ring of fractions}) of the ring $R$ (if it exists) and $Q^*$ is the group of units of $Q$;
        \item $\Den_l(R, \ga )$ is the set of left denominator sets $S$ of $R$ with $\ass (S)=\ga$ where $\ga$ is an ideal of $R$ and $\ass (S):= \{r\in R\, | \, sr=0$ for some $s\in S\}$,
             \item $\maxDen_l(R)$ is the set of maximal left denominator sets of $R$ (it is always a {\em non-empty} set, \cite{larglquot}).
\end{itemize}
{\bf Four new  criteria for a ring to have a semisimple left quotient ring}. Goldie's Theorem \cite{Goldie-PLMS-1960} characterizes left orders in semisimple rings, it is a criterion for a ring to have a semisimple  left quotient ring (earlier,  characterizations were given, by Goldie \cite{Goldie-PLMS-1958} and Lesieur and Croisot \cite{Lesieur-Croisot-1959}, of left orders in a simple Artinian ring).

\begin{theorem}\label{GoldieThm}
{\bf (Goldie's Theorem, \cite{Goldie-PLMS-1960})} A ring has a semisimple left quotient ring iff it is a semiprime  ring that satisfies the ascending chain condition on left annihilators and does not contain infinite direct sums of nonzero left ideals.
\end{theorem}

 In  \cite{larglquot} and  \cite{Bav-intdifline},  several new concepts
  are introduced (and studied): {\em the largest left quotient ring of a ring, the largest regular left Ore set of a ring, a  maximal left denominator set of a ring, the left localization radical of a ring, a left localization maximal ring} (see Section \ref{PRLM} for details). Their universal nature naturally leads to the present criteria for a ring to have a semisimple left quotient ring. In the paper, several new concepts are introduced that are used in proofs: {\em the core of an Ore set, the sets of left localizable and left non-localizable elements, the set of completely left localizable elements of a ring}.  The First Criterion is given via the set $\CM := \maxDen_l(R)$.

 \begin{itemize}
\item ({\bf Theorem \ref{31Dec12}, The First Criterion})  {\em A ring $R$ has a semisimple left quotient ring $Q$ iff $\CM $ is a finite set, $\bigcap_{S\in \CM } \ass (S) =0$ and,  for each $S\in \CM$,  the ring $S^{-1}R$ is a simple left Artinian ring. In this case, $Q\simeq \prod_{S\in \CM} S^{-1}R$.  }
\end{itemize}

 The Second Criterion  is given via the minimal primes of $R$ and certain explicit  multiplicative sets associated with them. On the one hand, the Second Criterion  stands between Goldie's Theorem and the First Criterion in terms how it is formulated. On the other hand, it goes further then the First Criterion  in the sense that it describes explicitly the maximal left denominator sets and the left quotient ring  of a ring with a semisimple left quotient ring.

 \begin{itemize}
\item ({\bf Theorem \ref{1Jan13}, The Second Criterion})  {\em  Let $R$ be a ring. The following statements are equivalent. }
\begin{enumerate}
\item {\em The ring $R$ has a semisimple left quotient ring $Q$.}
\item
\begin{enumerate}
\item {\em The ring $R$ is a semiprime ring.}
\item {\em The set $\Min (R)$ of minimal  primes of  $R$ is a finite set.}
    \item {\em For each $\gp\in \Min(R)$, the set $S_\gp := \{ c\in R\, | \, c+\gp\in \CC_{R/\gp} \}$ is a left denominator set of the ring $R$ with} $\ass (S_\gp )=\gp$.
        \item {\em  For each $\gp\in \Min(R)$, the ring $S_\gp^{-1} R$ is a simple left Artinian ring. }
\end{enumerate}
\end{enumerate}
{\em If one of the two equivalent conditions holds then} $\maxDen_l(R)=\{ S_\gp \, | \, \gp \in \Min (R)\}$ {\em and }   $Q\simeq \prod_{\gp\in \Min (R)}S_\gp^{-1}R$.
\end{itemize}

The Third Criterion (Theorem \ref{20Jan13}) can be seen as a `weak'  version of Goldie's Theorem in the sense that the conditions are `weaker' than that of Goldie's Theorem. In applications, it could be `easier' to verify whether a ring satisfies the conditions of Theorem \ref{20Jan13} comparing with Goldie's Theorem as Theorem \ref{20Jan13} `reduces' Goldie's Theorem essentially to the prime case and reveals the `local' nature of Goldie's Theorem.

 \begin{itemize}
\item ({\bf Theorem \ref{20Jan13}, The Third Criterion})  {\em  Let $R$ be a ring. The following statements are equivalent. }
\begin{enumerate}
\item {\em The ring $R$ has a semisimple left quotient ring $Q$.}
\item {\em The ring $R$ is a semiprime ring with $|\Min (R)|<\infty$ and, for each $\gp \in \Min (R)$, the ring $R/ \gp$ is a left Goldie ring.}
\end{enumerate}
\end{itemize}
The condition $|\Min (R)|<\infty$ can be replaced by any of the  four equivalent conditions of Theorem \ref{2.2.15-[MR]}, e.g., `the ring $R$ has a.c.c. on annihilator ideals.'

As  far as applications are concerned, Theorem \ref{20Jan13} has a useful corollary.

 \begin{itemize}
\item ({\bf Theorem \ref{22Jan13}, The Fourth Criterion})  {\em  Let $R$ be a ring. The following statements are equivalent. }
\begin{enumerate}
\item {\em The ring $R$ has a semisimple left quotient ring $Q$.}
\item {\em  There are left denominator sets $S_1', \ldots , S_n'$ of the ring $R$ such that  the rings $R_i:=S_i^{-1}R$, $i=1, \ldots , n$, are simple left Artinian rings and the map
    $$ \s := \prod_{i=1}^n \s_i : R\ra \prod_{i=1}^n R_i, \;\; R\mapsto (\frac{r}{1}, \ldots , \frac{r}{1}), $$
    is an injection where} $\s_i : R\ra R_i$, $r\mapsto \frac{r}{1}$.
\end{enumerate}
\end{itemize}

 The paper is organized as follows. In Section \ref{PRLM}, necessary concepts are introduced and results are collected that are used in the profs of Theorem \ref{31Dec12},  Theorem \ref{1Jan13},  Theorem \ref{20Jan13} and Theorem \ref{22Jan13}.

In Section \ref{FIRCRMD}, the proof of Theorem \ref{31Dec12} is given.

In Section \ref{TSECMP}, the proof of Theorem \ref{1Jan13} is given.

In Section \ref{THIGLCR}, the proof of Theorem \ref{20Jan13} is given.

In Section \ref{TFCCLD}, the proof of Theorem \ref{22Jan13} is given.

In Section \ref{LDSDFPR}, it is shown that the set of maximal left denominators of  a finite direct product $\prod_{i=1}^n R_i$ of rings is a union of the sets of maximal left denominator
sets of the rings $R_i$ (Theorem \ref{c26Dec12}).

In Section \ref{CRLSLG}, a criterion (Theorem \ref{24Feb13}) is given for the factor ring $R/ \gll$ (where $\gll$ is the left localization radical of $R$) to have a semisimple left quotient ring. The criterion is given in terms of the ring $R$ rather than of $R/ \gll$ and is based on four criteria above, it is a long statement. Let us give a flavour.

 \begin{itemize}
\item ({\bf Theorem \ref{24Feb13}})  {\em The following statements are equivalent. }
\begin{enumerate}
\item {\em The ring $R/\gll $ has a semisimple left quotient ring $Q$.}
    \item
\begin{enumerate}
\item  $|\maxDen_l(R)|<\infty$.
\item  {\em For every $S\in \maxDen_l(R)$, $S^{-1}R$ is a simple left Artinian ring}.
\end{enumerate}
\end{enumerate}
\end{itemize}
So, Theorem \ref{24Feb13} characterizes precisely  the class of rings that  have only finitely many maximal left denominators sets and all the left localizations at them are simple left Artinian rings.

The proofs of Theorem \ref{31Dec12} and Theorem \ref{1Jan13} are based on completely different ideas from existing proofs of Goldie's Theorem and Goldie's Theorem is not used. The key idea of the proof of Theorem \ref{20Jan13} is to reduce it to Theorem \ref{1Jan13} and the prime case of Goldie's Theorem. The key idea of the proof of Theorem \ref{22Jan13}  is to reduce it to Theorem \ref{20Jan13}.

 $\noindent $

 {\bf Old and new criteria for a ring to have a left Artinian left quotient ring}.
  Goldie's Theorem gives an answer to the question: {\em When a ring has a semisimple left quotient ring?} The next natural question (that was posed in 50s) and which is a generalization of the previous one: {\em Give a criterion for a ring to have a left Artinian left quotient ring.}
  Small \cite{Small-ArtQuotRings-66, Small-CorrArtQuotRings-66},    Robson \cite{Robson-ArtQuotRings-67}, and latter  Tachikawa \cite{Tachikawa-AutQuotR-71} and Hajarnavis \cite{Hajarnavis-ThmSmall-72}  have given different criteria  for a ring to have a left Artinian left quotient ring. Recently, the author \cite{Bav-genGoldie} has given three more criteria in the spirit of the present paper. We should mention contribution to the old criteria  of Talintyre, Feller and Swokowski \cite{Feller-Swokowski-RefMax-61}, \cite{Talintyre-QuotRinMin-66}. Talintyre \cite{Talintyre-QuotRinMax-63} and Feller and Swokowski \cite{Feller-Swokowski-RefPROCams-61} have given conditions
which are sufficient for a left Noetherian ring to
have a left quotient ring. Further, for a left Noetherian ring which has
a left quotient ring, Talintyre \cite{Talintyre-QuotRinMin-66} has established necessary and sufficient
conditions for the left quotient ring to be left Artinian. In the proofs of all these criteria (old and new) Goldie's Theorem is used. Each of the criteria comprises  several conditions. The conditions in the criteria of Small, Robson and Hajarnavis are `strong' and are given in terms of the ring $R$ rather than of its factor ring $\bR = R/ \gn $. On the contrary,  the conditions of the criteria in \cite{Bav-genGoldie}  are `weak' and given in terms of the ring $\bR$ and a {\em finite} set of explicit $\bR$-modules (they are certain explicit  subfactors of the prime radical $\gn$ of the ring $R$).



\section{Preliminaries}\label{PRLM}

In this section, we collect necessary results that are used in the proofs of this paper. More results on localizations of rings (and some of the missed standard definitions) the reader can find in \cite{Jategaonkar-LocNRings}, \cite{Stenstrom-RingQuot} and \cite{MR}.  In this paper the following notation will remained fixed.

$\noindent $

{\bf Notation}:

\begin{itemize}
\item $\Ore_l(R):=\{ S\, | \, S$ is a left Ore set in $R\}$; \item
$\Den_l(R):=\{ S\, | \, S$ is a left denominator set in $R\}$;
\item $\Loc_l(R):= \{ S^{-1}R\, | \, S\in \Den_l(R)\}$; \item
$\Ass_l(R):= \{ \ass (S)\, | \, S\in \Den_l(R)\}$ where $\ass
(S):= \{ r\in R \, | \, sr=0$ for some $s=s(r)\in S\}$;
\item $S_\ga=S_\ga (R)=S_{l,\ga }(R)$
 is the {\em largest element} of the poset $(\Den_l(R, \ga ),
\subseteq )$ and $Q_\ga (R):=Q_{l,\ga }(R):=S_\ga^{-1} R$ is  the
{\em largest left quotient ring associated to} $\ga$, $S_\ga $
exists (Theorem 
 2.1, \cite{larglquot});
\item In particular, $S_0=S_0(R)=S_{l,0}(R)$ is the largest
element of the poset $(\Den_l(R, 0), \subseteq )$, i.e. the largest regular  left Ore set of $R$,  and
$Q_l(R):=S_0^{-1}R$ is the largest left quotient ring of $R$ \cite{larglquot}  ;
\item $\Loc_l(R, \ga ):= \{ S^{-1}R\, | \, S\in \Den_l(R, \ga
)\}$.
\end{itemize}

In \cite{larglquot}, we introduce the following new concepts and
prove their existence for an {\em arbitrary} ring: {\em the
largest left quotient ring of a ring, the largest  regular left
 Ore  set of a ring, the left localization radical of a ring, a maximal left denominator set,  a maximal left quotient ring of a ring,
 a (left) localization
maximal ring}.
Using an analogy with rings, the counter parts of these concepts for rings would be a left maximal ideal, the Jacobson radical, a simple factor ring. These concepts turned out to be very useful in Localization Theory and Ring Theory. They allowed us to look at old/classical results from  a new more general perspective and to give new equivalent statements  to the classical results using a new language and a new approach as the present paper, \cite{larglquot}, \cite{Bav-intdifline}, \cite{Bav-genGoldie},  \cite{Bav-LocArtRing} and  \cite{Bav-locmaxrings}  and several other papers under preparation demonstrate.

{\bf The largest regular left Ore set and the largest left
quotient ring of a ring}. Let $R$ be a ring. A {\em
multiplicatively closed subset} $S$ of $R$ or a {\em
 multiplicative subset} of $R$ (i.e. a multiplicative sub-semigroup of $(R,
\cdot )$ such that $1\in S$ and $0\not\in S$) is said to be a {\em
left Ore set} if it satisfies the {\em left Ore condition}: for
each $r\in R$ and
 $s\in S$, $ Sr\bigcap Rs\neq \emptyset $.
Let $\Ore_l(R)$ be the set of all left Ore sets of $R$.
  For  $S\in \Ore_l(R)$, $\ass (S) :=\{ r\in
R\, | \, sr=0 \;\; {\rm for\;  some}\;\; s\in S\}$  is an ideal of
the ring $R$.


A left Ore set $S$ is called a {\em left denominator set} of the
ring $R$ if $rs=0$ for some elements $ r\in R$ and $s\in S$ implies
$tr=0$ for some element $t\in S$, i.e. $r\in \ass (S)$. Let
$\Den_l(R)$ be the set of all left denominator sets of $R$. For
$S\in \Den_l(R)$, let $S^{-1}R=\{ s^{-1}r\, | \, s\in S, r\in R\}$
be the {\em left localization} of the ring $R$ at $S$ (the {\em
left quotient ring} of $R$ at $S$). Let us stress that in Ore's method of localization one can localize {\em precisely} at left denominator sets.

In general, the set $\CC$ of regular elements of a ring $R$ is
neither left nor right Ore set of the ring $R$ and as a
 result neither left nor right classical  quotient ring ($Q_{l,cl}(R):=\CC^{-1}R$ and
 $Q_{r,cl}(R):=R\CC^{-1}$) exists.
 Remarkably, there  exists the largest
 regular left Ore set $S_0= S_{l,0} = S_{l,0}(R)$, \cite{larglquot}. This means that the set $S_{l,0}(R)$ is an Ore set of
 the ring $R$ that consists
 of regular elements (i.e., $S_{l,0}(R)\subseteq \CC$) and contains all the left Ore sets in $R$ that consist of
 regular elements. Also, there exists the largest regular (left and right) Ore set $S_{l,r,0}(R)$ of the ring $R$.
 In general, all the sets $\CC$, $S_{l,0}(R)$, $S_{r,0}(R)$ and $S_{l,r,0}(R)$ are distinct, for example,
 when $R= \mI_1= K\langle x, \der , \int\rangle$  is the ring of polynomial integro-differential operators  over a field $K$ of characteristic zero,  \cite{Bav-intdifline}. In  \cite{Bav-intdifline},  these four sets are found for $R=\mI_1$.

$\noindent $

{\it Definition}, \cite{Bav-intdifline}, \cite{larglquot}.    The ring
$$Q_l(R):= S_{l,0}(R)^{-1}R$$ (respectively, $Q_r(R):=RS_{r,0}(R)^{-1}$ and
$Q(R):= S_{l,r,0}(R)^{-1}R\simeq RS_{l,r,0}(R)^{-1}$) is  called
the {\em largest left} (respectively, {\em right and two-sided})
{\em quotient ring} of the ring $R$.

$\noindent $

 In general, the rings $Q_l(R)$, $Q_r(R)$ and $Q(R)$
are not isomorphic, for example, when $R= \mI_1$, \cite{Bav-intdifline}.  The next
theorem gives various properties of the ring $Q_l(R)$. In
particular, it describes its group of units.


\begin{theorem}\label{4Jul10}
\cite{larglquot}
\begin{enumerate}
\item $ S_0 (Q_l(R))= Q_l(R)^*$ {\em and} $S_0(Q_l(R))\cap R=
S_0(R)$.
 \item $Q_l(R)^*= \langle S_0(R), S_0(R)^{-1}\rangle$, {\em i.e. the
 group of units of the ring $Q_l(R)$ is generated by the sets
 $S_0(R)$ and} $S_0(R)^{-1}:= \{ s^{-1} \, | \, s\in S_0(R)\}$.
 \item $Q_l(R)^* = \{ s^{-1}t\, | \, s,t\in S_0(R)\}$.
 \item $Q_l(Q_l(R))=Q_l(R)$.
\end{enumerate}
\end{theorem}

{\bf The maximal denominator sets and the maximal left localizations  of a ring}. The set $(\Den_l(R), \subseteq )$ is a poset (partially ordered
set). In \cite{larglquot}, it is proved  that the set
$\maxDen_l(R)$ of its maximal elements is a {\em non-empty} set.

$\noindent $

{\it Definition}, \cite{larglquot}. An element $S$ of the set
$\maxDen_l(R)$ is called a {\em maximal left denominator set} of
the ring $R$ and the ring $S^{-1}R$ is called a {\em maximal left
quotient ring} of the ring $R$ or a {\em maximal left localization
ring} of the ring $R$. The intersection
\begin{equation}\label{llradR}
\gll_R:=\llrad (R) := \bigcap_{S\in \maxDen_l(R)} \ass (S)
\end{equation}
is called the {\em left localization radical } of the ring $R$,
\cite{larglquot}.

$\noindent $

 For a ring $R$, there is the canonical exact
sequence 
\begin{equation}\label{llRseq}
0\ra \gll_R \ra R\stackrel{\s }{\ra} \prod_{S\in \maxDen_l(R)}S^{-1}R, \;\; \s := \prod_{S\in \maxDen_l(R)}\, \s_S,
\end{equation}
where $\s_S:R\ra S^{-1}R$, $r\mapsto \frac{r}{1}$. For a ring $R$ with a semisimple left quotient ring, Theorem \ref{1Jan13} shows that the left localization radical $\gll_R$ coincides with the prime radical $\gn_R$ of $R$: $\gll_R=\bigcap_{\gp\in \Min (R)}\gp =\gn_R= 0$. In general, this  is not the case even for left Artinian rings \cite{Bav-LocArtRing}.

$\noindent $

{\it Definition}. The sets
$$ \CL_l(R):= \bigcup_{S\in \maxDen_l(R)}S\;\; {\rm and}\;\; \CN\CL_l(R):=R\backslash \CL_l(R)$$ are called the sets of {\em left localizable} and {\em left non-localizable elements} of $R$, respectively, and the intersection
$$\CC_l(R):=\bigcap_{S\in \maxDen_l(R)} S$$
 is called the {\em set of completely left localizable elements} of $R$. Clearly, $R^*\subseteq \CC_l(R)$.

$\noindent $

{\bf The maximal elements of $\Ass_l(R)$}.  Let $\maxAss_l(R)$ be
the set of maximal elements of the poset $(\Ass_l(R), \subseteq )$
and

\begin{equation}\label{mADen}
\assmaxDen_l(R) := \{ \ass (S) \, | \, S\in \maxDen_l(R) \}.
\end{equation}
These two sets are equal (Proposition \ref{b27Nov12}), a proof is
based on Lemma \ref{1a27Nov12} and Corollary \ref{d4Jan13}. Recall that for an non-empty set $X$ or $R$, $\lann (X):=\{ r\in R\, | \, rX=0\}$ is the {\em left annihilator} of the set $X$, it is a left ideal of $R$.

\begin{lemma}\label{1a27Nov12}
 \cite{larglquot}
Let $S\in \Den_l(R, \ga )$ and $T\in \Den_l(R, \gb )$ be such that $ \ga \subseteq \gb$. Let $ST$ be the multiplicative semigroup generated by $S$ and $T$ in $(R,\cdot )$.  Then
\begin{enumerate}
\item $\lann (ST)\subseteq \gb$.
\item $ST \in \Den_l(R, \gc )$ and $\gb \subseteq \gc$.
\end{enumerate}
\end{lemma}

\begin{corollary}\label{d4Jan13}
Let $R$ be a ring, $S\in \maxDen_l(R)$ and $T\in \Den_l(R)$. Then
\begin{enumerate}
\item $T\subseteq S$ iff $\ass (T)\subseteq \ass (S)$.
\item If, in addition, $T\in \maxDen_l(R)$ then $S=T$ iff $\ass (S) = \ass (T)$.
\end{enumerate}
\end{corollary}

{\it Proof}. 1.  $(\Rightarrow )$ If $T\subseteq S$ then $\ass (T)\subseteq \ass (S)$.

$(\Leftarrow )$ If $\ass (T)\subseteq \ass (S)$. then, by Lemma \ref{1a27Nov12}, $ST\in \Den_l(R)$ and $S\subseteq ST$, hence $S= ST$, by the maximality of $S$. Then $T\subseteq S$.

2. Statement 2 follows from statement 1.  $\Box $


\begin{proposition}\label{b27Nov12}
\cite{larglquot} $\; \maxAss_l(R)= \assmaxDen_l(R)\neq \emptyset$. In particular, the ideals of this set are incomparable (i.e. neither $\ga\nsubseteq \gb$ nor $\ga\nsupseteq \gb$).
\end{proposition}

{\bf Properties of the maximal left quotient rings of a ring}.
The next theorem describes various properties of the maximal left
quotient rings of a ring. In particular, their groups of units and
their largest left quotient rings. The key moment in the proof is to use Theorem \ref{4Jul10}.

\begin{theorem}\label{15Nov10}
\cite{larglquot} Let $S\in \maxDen_l(R)$, $A= S^{-1}R$, $A^*$ be
the group of units of the ring $A$; $\ga := \ass (S)$, $\pi_\ga
:R\ra R/ \ga $, $ a\mapsto a+\ga$, and $\s_\ga : R\ra A$, $
r\mapsto \frac{r}{1}$. Then
\begin{enumerate}
\item $S=S_\ga (R)$, $S= \pi_\ga^{-1} (S_0(R/\ga ))$, $ \pi_\ga
(S) = S_0(R/ \ga )$ and $A= S_0( R/\ga )^{-1} R/ \ga = Q_l(R/ \ga
)$. \item  $S_0(A) = A^*$ and $S_0(A) \cap (R/ \ga )= S_0( R/ \ga
)$. \item $S= \s_\ga^{-1}(A^*)$. \item $A^* = \langle \pi_\ga (S)
, \pi_\ga (S)^{-1} \rangle$, i.e. the group of units of the ring
$A$ is generated by the sets $\pi_\ga (S)$ and $\pi_\ga^{-1}(S):=
\{ \pi_\ga (s)^{-1} \, | \, s\in S\}$. \item $A^* = \{ \pi_\ga
(s)^{-1}\pi_\ga ( t) \, |\, s, t\in S\}$. \item $Q_l(A) = A$ and
$\Ass_l(A) = \{ 0\}$.     In particular, if $T\in \Den_l(A, 0)$
then  $T\subseteq A^*$.
\end{enumerate}
\end{theorem}

{\bf The left localization maximal rings}. These are precisely the rings  in which  we cannot invert anything on the left (in the sense of Ore).

 $\noindent $

 {\it Definition}, \cite{larglquot}. A ring $A$ is
called a {\em left localization maximal ring} if $A= Q_l(A)$ and
$\Ass_l(A) = \{ 0\}$. A ring $A$ is called a {\em right
localization maximal ring} if $A= Q_r(A)$ and $\Ass_r(A) = \{
0\}$. A ring $A$ which is a left and right localization maximal
ring is called a {\em (left and right) localization maximal ring}
(i.e. $Q_l(A) =A=Q_r(A)$ and $\Ass_l(A) =\Ass_r(A) = \{ 0\}$).

$\noindent $


{\it Example}. Let $A$ be a simple ring. Then $Q_l(A)$ is a left
localization maximal  ring and $Q_r(A)$ is a right localization
maximal ring. In particular, a division ring is a (left and right) localization maximal ring. More generally, a simple left Artinian ring (i.e. the matrix ring over a division ring) is a (left and right)
localization maximal ring.

$\noindent $

Let $\maxLoc_l(R)$ be the set of maximal elements of the poset
$(\Loc_l(R), \ra )$ where $A\ra B$ for $A,B\in \Loc_l(R)$ means that there exist $S,T\in \Den_l(R)$ such that $S\subseteq T$, $A\simeq S^{-1}R$  and $B\simeq T^{-1}R$ (then there exists a natural ring homomorphism $A\ra B$, $s^{-1}r\mapsto s^{-1}r$). Then (see \cite{larglquot}),
\begin{equation}\label{mADen1}\maxLoc_l(R) = \{ S^{-1}R \, | \, S\in \maxDen_l(R) \}= \{ Q_l(R/
\ga ) \, | \, \ga \in \assmaxDen_l(R)\}.
\end{equation}

The next theorem is a criterion of when  a left quotient ring of a
ring is a maximal left quotient ring.

\begin{theorem}\label{21Nov10}
\cite{larglquot} Let  a ring $A$ be a left localization of a ring
$R$, i.e. $A\in \Loc_l(R, \ga )$ for some $\ga \in \Ass_l( R)$.
Then $A\in \maxLoc_l( R)$ iff $Q_l( A) = A$ and  $\Ass_l(A) = \{
0\}$, i.e. $A$ is a left localization maximal ring.
\end{theorem}


Theorem \ref{21Nov10} shows that the left localization maximal
rings are precisely the localizations of all the rings at their
maximal left denominators sets.

$\noindent $

{\bf The core of a left Ore set}. The following definition is one of the key concepts that is used in the proof of the First Criterion (Theorem \ref{31Dec12}).

$\noindent $

{\it Definition}, \cite{Bav-locmaxrings}. Let $R$ be a ring and $S\in \Ore_l(R)$. The {\em core} $S_c$ of the left Ore set $S$ is the set of all the elements $s\in S$ such that $\ker (s\cdot) = \ass (S)$ where $s\cdot : R\ra R$, $r\mapsto sr$.  


\begin{lemma}\label{b2Dec12}
If $S\in \Den_l(R)$ and $S_c\neq \emptyset$ then
\begin{enumerate}
\item $SS_c\subseteq S_c$.
\item For any $s\in S$ there exists an element $t\in S$ such that $ts \in S_c$.
\end{enumerate}
\end{lemma}

{\it Proof}. 1. Trivial.

2. Statement 2 follows directly from the left Ore condition: fix
an element $s_c\in S_c$, then $ts= rs_c$ for some elements $t\in
S$ and $r\in R$. Since $\ass (S)\supseteq \ker ( ts\cdot
)=\ker(rs_c\cdot ) \supseteq \ker (s_c\cdot ) = \ass (S)$, i.e.
$\ker (ts\cdot ) = \ass (S)$, we have $ts \in S_c$.  $\Box $


\begin{theorem}\label{A2Dec12}
Suppose that $S\in \Den_l(R, \ga )$ and $S_c\neq \emptyset$. Then
\begin{enumerate}
\item $S_c\in \Den_l(R, \ga )$.
\item The map $\th : S_c^{-1}R\ra S^{-1}R$, $s^{-1}r\mapsto s^{-1}r$, is a ring isomorphism. So, $S_c^{-1}R\simeq S^{-1}R$.
\end{enumerate}
\end{theorem}

{\it Proof}. 1. By Lemma \ref{b2Dec12}.(1), $S_cS_c\subseteq S_c$,
 and so  the set $S_c$ is a multiplicative set. By Lemma
\ref{b2Dec12}.(2), $S_c\in \Ore_l(R)$: for any elements $s_c\in
S_c$ and $r\in R$, there are elements $s\in S$ and $r'\in R$ such
that $sr = r's_c$ (since $S\in \Ore_l(R)$). By Lemma
\ref{b2Dec12}.(2), $s_c':= ts\in S_c$ for some $t\in S$. Then $s_c' r=
tr' s_c$.

If $rs_c=0$ for some elements $r\in R$ and $s_c\in S_c$ then $sr=0$ for some element $s\in S$ (since $S\in \Ore_l(R)$). By Lemma \ref{b2Dec12}.(2), $s_c':=ts\in S_c$ for some element $t\in S$, hence $s_c'r=0$. Therefore, $S_c\in \Den_l(R, \ga )$.

2. By statement 1 and the universal property of left Ore
localization,  the map $\th$ is a well-defined monomorphism. By
Lemma \ref{b2Dec12}.(2), $\th$ is also a surjection: let
$s^{-1}r\in S^{-1}R$, and $s_c:= ts\in S_c$ for some element $t\in
S$ (Lemma \ref{b2Dec12}.(2)). Then
$$ s^{-1}r= s^{-1}t^{-1}tr= (ts)^{-1}tr= s_c^{-1} tr. \;\; \Box $$


The core of every maximal left denominator set of a semiprime left Goldie ring is a non-empty set and Theorem \ref{31Dec12}.(6) gives its  explicit description (via the minimal primes).

$\noindent $

{\bf The maximal left quotient rings of a finite direct product of rings}.
\begin{theorem}\label{c26Dec12}
 Let $R=\prod_{i=1}^n R_i$ be a direct product of rings $R_i$. Then
for each $i=1, \ldots , n$, the map
\begin{equation}\label{1aab1}
\maxDen_l(R_i) \ra \maxDen_l(R), \;\; S_i\mapsto R_1\times\cdots \times S_i\times\cdots \times R_n,
\end{equation}
is an injection. Moreover, $\maxDen_l(R)=\coprod_{i=1}^n \maxDen_l(R_i)$ in the sense of (\ref{1aab1}), i.e.
$$ \maxDen_l(R)=\{ S_i\, | \, S_i\in \maxDen_l(R_i), \; i=1, \ldots , n\},$$
$S_i^{-1}R\simeq S_i^{-1}R_i$, $\ass_R(S_i)= R_1\times \cdots \times \ass_{R_i}(S_i)\times\cdots \times R_n$. The core of the left denominator set $S_i$ in $R$ coincides with the core $S_{i,c}$ of the left denominator set $S_i$ in $R_i$, i.e.
$$(R_1\times\cdots \times S_i\times\cdots \times R_n)_c=0\times\cdots \times S_{i,c}\times\cdots \times 0.$$
\end{theorem}

The proof of Theorem \ref{c26Dec12} is given in Section \ref{LDSDFPR}.

$\noindent $

{\bf A bijection between $\maxDen_l(R)$ and $\maxDen_l(Q_l(R))$}.
\begin{proposition}\label{A8Dec12}
  Let $R$ be a ring, $S_l$ be the  largest regular left Ore set of the ring $R$, $Q_l:= S_l^{-1}R$ be the largest left quotient ring of the ring $R$, and $\CC$ be the set of regular elements of the ring $R$. Then
\begin{enumerate}
\item $S_l\subseteq S$ for all $S\in \maxDen_l(R)$. In particular,
$\CC\subseteq S$ for all $S\in  \maxDen_l(R)$ provided $\CC$ is a
left Ore set. \item Either $\maxDen_l(R) = \{ \CC \}$ or,
otherwise, $\CC\not\in\maxDen_l(R)$. \item The map $$
\maxDen_l(R)\ra \maxDen_l(Q_l), \;\; S\mapsto SQ_l^*=\{ c^{-1}s\,
| \, c\in S_l, s\in S\},
$$ is a bijection with the inverse $\CT \mapsto \s^{-1} (\CT )$
where $\s : R\ra Q_l$, $r\mapsto \frac{r}{1}$, and $SQ_l^*$ is the
sub-semigroup of $(Q_l, \cdot )$ generated by the set  $S$ and the
group $Q_l^*$ of units of the ring $Q_l$, and $S^{-1}R= (SQ_l^*)^{-1}Q_l$.
    \item  If $\CC$ is a left Ore set then the map $$ \maxDen_l(R)\ra \maxDen_l(Q), \;\; S\mapsto SQ^*=\{ c^{-1}s\,
| \, c\in \CC, s\in S\}, $$ is a bijection with the inverse $\CT
\mapsto \s^{-1} (\CT )$ where $\s : R\ra Q$, $r\mapsto
\frac{r}{1}$, and $SQ^*$ is the sub-semigroup of $(Q, \cdot )$
generated by the set  $S$ and the group $Q^*$ of units of the ring
$Q$, and $S^{-1}R= (SQ^*)^{-1}Q$.
\end{enumerate}
\end{proposition}

{\it Proof}. 1. Let $S\in \maxDen_l(R)$. By Lemma \ref{1a27Nov12},
$S_l\subseteq S$.

2. Statement 2 follows from statement 1.

3. Statement 3 follows  from statement 1 and Proposition 3.4.(1),
\cite{larglquot}.

4. If $\CC\in \Ore_l(R)$ then $\CC = S_l$ and statement 4 is a particular case of statement 3.  $\Box $



\section{The First Criterion (via the maximal left denominator sets)}\label{FIRCRMD}

The aim of this section is to give a criterion for a ring $R$ to have  a semisimple left quotient ring (Theorem \ref{31Dec12}). The implication $(\Leftarrow )$ is the most difficult  part of Theorem \ref{31Dec12}. Roughly speaking, it proceeds by establishing
 properties 1-9 which are interesting on their own  and constitute the structure of the proof. These properties show that the relationships between the ring $R$ and its semisimple left quotient ring  $Q$  are as natural as possible and are as simple as possible (`simple' in the sense that the connections between the properties of $R$ and $Q$ are strong).

\begin{theorem}\label{31Dec12}
({\bf The First Criterion}) A ring $R$  have a semisimple left quotient ring $Q$ iff $\maxDen_l(R)=\{ S_1, \ldots , S_n\}$ is a finite set, $\bigcap_{i=1}^n \ass (S_i)=0$ and $R_i:=S_i^{-1}R$ is a simple left Artinian ring for $i=1, \ldots , n$. If one of the equivalent conditions hold then
the map $$ \s := \prod_{i=1}^n \s_i : R\mapsto Q':=\prod_{i=1}^n R_i, \;\; r\mapsto (r_1, \ldots , r_n), $$ is a ring monomorphism where $r_i=\s_i(r)$ and $\s_i: R\ra R_i$, $r\mapsto \frac{r}{1}$;  and
\begin{enumerate}
\item $\CC = \bigcap_{i=1}^nS_i$.
\item The map $\s':= \prod_{i=1}^n \s_i': Q\ra Q'$, $ c^{-1}r\mapsto (c^{-1}r, \ldots , c^{-1}r)$, is a ring isomorphism where $\s_i':Q\ra Q'$, $c^{-1}r\mapsto c^{-1}r$. We identify the rings $Q$ and $Q'$ via $\s'$.
\item $\maxDen_l(Q')=\{ S_1', \ldots , S_n'\}$  where $S_i':=R_1\times\cdots \times R_i^*\times \cdots \times R_n$, $R_i^*$ is the group of units of the ring $R_i$, $\ga_i':= \ass (S_i') = R_1\times\cdots \times 0\times \cdots \times R_n$, $S_i'^{-1}Q'\simeq R_i$ for $i=1, \ldots , n$. The core $S_{i, c}'$ of the left denominator set  $S_i'$ is equal to $R_i^*=0\times \cdots \times 0 \times R_i^*\times 0 \times \cdots \times 0$.
    \item The map $\maxDen_l(Q')\ra \maxDen_l(R)$, $S_i'\mapsto S_i:= R\cap S_i'= \s_i^{-1}(R_i^*)$ is a bijection, $\ga_i:= \ass (S_i) = R\cap \ga_i'$.
        \item For all $i=1, \ldots , n$, $S_i\nsubseteq \bigcup_{j\neq i} S_j$. Moreover,
        $\emptyset \neq S_{i,c}\subseteq S_i\backslash \bigcup_{j\neq i} S_j= R\cap (R_i^*\times \prod_{j\neq i} R_i^0)$ where $R_i^0:= R_i\backslash R_i^*$ is the set of zero divisors of the ring $R_i$.
        \item $S_{i,c}=R\cap S_{i,c}'=S_i\cap S_{i,c}'=S_i\cap \bigcap_{j\neq i} \ga_j= R_i^*\cap \bigcap_{j\neq i} \ga_j= (\bigcap_{j\neq i} \ga_j) \backslash R_i^0\neq \emptyset $ for $i=1, \ldots , n$ where $R_i^*, R_i^0\subseteq \prod_{j=1}^n  R_j$ are  the natural inclusions $r\mapsto (0, \ldots , 0, r,0, \ldots , 0)$. For all $i\neq j$, $S_{i,c}S_{j,c}=0$.
            \item $\CC':= S_{1, c}+\cdots +S_{n,c}\in \Den_l(R, 0)$, $\CC'^{-1}R\simeq Q$, $Q=\{\sum_{i=1}^n s_i^{-1}a_i\, | \, s_i\in S_{i, c}, a_i\in \bigcap_{j\neq i} \ga_j$ for $i=1, \ldots , n\}$, $\CC\CC'\subseteq \CC'$ and $\CC S_{i,c}\subseteq S_{i,c}$ for  $i=1, \ldots , n$.
                \item $Q^*=\{ s^{-1}t\, | \, s,t\in \CC'\}= \{ s^{-1}t\, | \, s,t\in \CC \}$.
                    \item $\CC= \{ s^{-1}t\, | \, s,t\in \CC', s^{-1}t\in R\}$.
\end{enumerate}
\end{theorem}

{\it Proof}.  $(\Rightarrow )$  Suppose that the left quotient ring $Q$ of the ring $R$ is a semisimple ring, i.e.
$$Q\simeq \prod_{i=1}^n R_i$$ where $R_i$ are simple left Artinian rings. Every simple left Artinian ring is a left localization maximal ring, hence $\maxDen_l(R_i) = \{R_i^*\}$. Then, by Theorem \ref{c26Dec12},  $$\maxDen_l(Q) =\maxDen_l(\prod_{i=1}^nR_i)=\{ S_1', \ldots , S_n'\}$$ where $S_i'=R_1\times\cdots \times R_i^*\times \cdots \times R_n$, $$\ass_{Q}(S_i')=R_1\times \cdots \times 0 \times \cdots \times R_i\;\; {\rm  and}\;\;  S_i'^{-1}Q\simeq Q/\ass_{Q}(S_i')\simeq R_i.$$ The core $S_{i,c}'$ of the maximal  left denominator set $S_i'$ is  $R_i^*=0\times \cdots \times 0 \times R_i^*\times 0 \times \cdots \times 0$.
  The ring $R\ra Q$, $r\mapsto \frac{r}{1}$, is a ring monomorphism. We identify the ring $R$ with its image in $Q$. By Proposition \ref{A8Dec12}.(4), the map
 $$ \maxDen_l(Q)\ra \maxDen_l(R), \;\; S_i'\mapsto S_i:= R\cap S_i',$$ is a bijection and  $S_i^{-1}R\simeq S_i'^{-1}Q\simeq R_i$.  The inclusions  $\ass_{R}(S_i)\subseteq \ass_{Q}(S_i')$ where $i=1, \ldots , n$ imply   $\bigcap_{i=1}^n\ass_{R}(S_i)\subseteq \bigcap_{i=1}^n \ass_{Q}(S_i')=0$, and so $\bigcap_{i=1}^n\ass_{R}(S_i)=0$.

$(\Leftarrow )$  Suppose that $\maxDen_l(R)=\{S_1, \ldots , S_n\}$, $\bigcap_{i=1}^n\ass (S_i) =0$ and $R_i:=S_i^{-1}R$ is a simple left Artinian ring for $i=1, \ldots , n$. We keep the notation of the theorem. Then $\s$ is a monomorphism (since  $\bigcap_{i=1}^n \ga_i =0$ where $\ga_i:=\ass (S_i)$), and we can identify the ring $R$ with its image in the ring $Q'$. Repeating word for word the arguments of the proof of the implication  $(\Rightarrow )$, by simply replacing the letter $Q$ by $Q'$,  we see that statement 3 holds (clearly, $S_{i,c}'= 0\times \cdots \times 0 \times R_i^*\times 0 \times \cdots \times 0$).  The group $Q'^*$ of units of the direct product $Q'$ is equal to $\prod_{i=1}^n R_i^*$.

{\em Step 1}: $\CC = R\bigcap Q'^*$: The inclusion $\CC \supseteq R\bigcap Q'^*$ is obvious. To prove that the reverse inclusion holds it suffices to show that, for each element $c\in \CC$, the map $$\cdot c:Q'\ra Q', \;\; q\mapsto qc,$$ is an injection (then necessarily, the map $\cdot c$ is a bijection, then  it is an automorphism of the left Artinian $Q'$-module $Q'$, its inverse is also an element of the type $\cdot c'$ for some element $c'\in Q'^*$. Then $c=(c')^{-1}\in Q'^*$, as required).  Suppose that $\ker_{Q'}(\cdot c)\neq 0$, we seek a contradiction. Fix a nonzero element $q=(s_1^{-1}r_1, \ldots , s_n^{-1}r_n)\in \ker_{Q'}(\cdot c)$. Without loss of generality we may assume that $s_1=\cdots = s_n=1$, multiplying several times, if necessary,
by  well-chosen elements  of the set $\bigcup_{i=1}^nS_i$ in order to get rid of the denominators. There is a nonzero component  of the element $q=(r_1, \ldots , r_n)$, say $r_i$. Then $r_i=\s_i(r)$ for some element $r\not\in \ga_i$ (otherwise, we would have $r_i=0$). The equality $ qc=0$ implies that $0=r_i\s_i(c)= \s_i(rc)$. Then $s_ir_ic=0$ for some element $s_i\in S_i$, and so $s_ir=0$ since $c\in \CC$. This means that $ 0=\s_i(r) = r_i$, a contradiction. The proof of Step 1 is complete.

{\em Step 2}: {\em For all} $i\neq j$, $S_i^{-1}\ga_j=R_j$: The ring $R_i = S_i^{-1}R$ is a left Artinian ring and $\ga_j$ is an ideal of the ring $R$, hence $S_i^{-1}\ga_j$ is an ideal of the simple ring $R_i$. There are two options: either $S_i^{-1}\ga_j =0$ or, otherwise,  $S_i^{-1}\ga_j =R_i$. Suppose that $S_i^{-1}\ga_j =0$, i.e. $\ga_j\subseteq \ga_i$, but this is not possible, by Proposition \ref{b27Nov12}.  Therefore, $S_i^{-1}\ga_j =R_i$.

{\em Step 3}: {\em For all} $i=1, \ldots , n$, $S_i\cap \bigcap_{j\neq i} \ga_j\neq \emptyset$:  For each pair of distinct indices  $i\neq j$, $R_i =S_i^{-1}\ga_j$ (by Step 2), and so $R_i\ni 1=s_i^{-1} a_j$ for some elements $s_i\in S_i$ and $a_j\in \ga_j$. Then $s_i-a_j\in \ga_i$, and so  there is an element $s_{ij}\in S_i$ such that $s_{ij}(s_i-a_j)=0$, and so
$$t_{ij}:= s_{ij}s_i = s_{ij}a_j\in S_i\cap \ga_j.$$
Then
\begin{equation}\label{tiSiaj}
t_i:=\prod_{j\neq i}t_{ij}\in S_i\cap \bigcap_{j\neq i}\ga_j.
\end{equation}

4. We have the commutative diagram of ring homomorphisms:

\begin{equation}\label{nRQsRi}
\xymatrix{
R\ar[rd]^{\s_i}\ar[r]^{\s } & Q'\ar[d]^{p_i}\\
 &  R_i}
\end{equation}

$\s_i = p_i\s$ where $p_i:Q'=\prod_{j=1}^n R_j\ra R_i$, $(r_1, \ldots , r_n) \mapsto r_i$.  By statement 3, $S_i^{-1}R=R_i\simeq S_i'^{-1}Q'$ where  $S_i\in \maxDen_l(R)$ and  $S_i'\in \maxDen_l(Q')$. Then, by Theorem \ref{15Nov10}.(3),  $S_i = \s_i^{-1}(R_i^*)$ and $ S_i'=p_i^{-1}(R_i^*)$. Therefore,
$$ S_i= \s_i^{-1} (R_i^*) = (p_i\s)^{-1} (R_i^*) = \s^{-1}(p_i^{-1}(R_i^*)) =\s^{-1}(S_i'),$$ and so the map
$$ \maxDen_l(Q')\ra \maxDen_l(R), \;\; S_i'\ra S_i = \s^{-1}(S_i') = R\cap S_i',$$
is a bijection. It follows from the commutative diagram (\ref{nRQsRi})  and Step 3 that
\begin{equation}\label{nRQsRi1}
 \ga_i = \ass (S_i) = R\cap (R_1\times \cdots \times 0 \times \cdots \times R_n) = R\cap \ga_i'.
\end{equation}
 In more detail, $\ga_i\subseteq R\cap \ga_i'$, by (\ref{nRQsRi}). Then $t_i\in S_i\cap \bigcap_{j\neq i} \ga_j'$, by (\ref{tiSiaj}), and so $t_i (R\cap \ga_i')=0$ (since $t_i(R\cap \ga_i')\subseteq \cap_{i=1}^n \ga_i'=0$). This equality implies the inclusion $ R\cap \ga_i'\subseteq \ga_i$, and so (\ref{nRQsRi1}) holds.

6. We claim that $S_{i,c} = R\cap S_{i,c}'$. The inclusion $S_{i,c} \supseteq  R\cap S_{i,c}'$ follows from (\ref{nRQsRi1}). Suppose that $S_{i,c}\backslash R\cap S_{i,c}'\neq \emptyset$, we seek a contradiction. Fix an element $s\in S_{i,c}\backslash R\cap S_{i,c}'$. Then $s=(s_1, \ldots ,s_n)$ with $s_j\neq 0$ for some $j$ such that $j\neq i$. Let $t_j$ be the element from Step 3, see (\ref{tiSiaj}),  i.e.
$$t_j\in S_j\cap \bigcap_{k\neq j}\ga_k\subseteq \ga_i,$$
 since $j\neq i$. Notice that $t_j=(0, \ldots , 0,t_j,0,\ldots ,0)$ and $t_j\in S_j\subseteq S_j'$, and so $t_j\in R_j^*$. On the  one hand, $st_j=0$ since $s\in S_i$ and $t_j\in \ga_i$, and so $0=\s_j(st_j) = s_jt_j$. On the other hand, $s_jt_j\neq 0$ since $s_j\neq 0$ and $t_j\in R_j^*$, a contradiction.  Therefore, $S_{i,c} = R\cap S_{i,c}'$. Then $S_{i,c}S_{j,c}=0$ for all $i\neq j$.  Then, for all $i=1, \ldots , n$,
\begin{eqnarray*}
 S_{i,c}&=& R\cap S_{i,c}'= S_i\cap S_{i,c}'=R\cap S_i'\cap \bigcap_{j\neq i}\ga_j'=(R\cap S_i')\cap \bigcap_{j\neq i} R\cap \ga_j'= S_i\cap\bigcap_{j\neq i} \ga_j \ni t_i,\\
  S_{i,c}&=& S_{i,c}'\cap \bigcap_{j\neq i} \ga_j =R_i^*\cap \bigcap_{j\neq i} \ga_j = (\bigcap_{j\neq i} \ga_j)\backslash R_i^0\neq \emptyset .\\
\end{eqnarray*}
So, statement 6 has been proven.

2 {\em and} 7. By statement 6, $\CC'\subseteq Q'^*=\prod_{i=1}^nR_i^*$, and so
$$\CC'\subseteq R\cap Q'^*= \CC ,$$
 by Step 1.  Since $\CC\subseteq Q'^*$ (Step 1), for each $i=1, \ldots , n$, $\CC S_{i,c}\subseteq S_{i,c}$, by statement 6. Then $\CC \CC'\subseteq \CC'$ since $\CC'=\sum_{i=1}^n S_{i,c}$.

By Theorem \ref{A2Dec12}.(2), $$S_{i,c}^{-1}R= S_i^{-1}R=R_i$$ for $i=1, \ldots , n$. Each element $q\in Q'$ can be written  as $q=(s_1^{-1}r_1, \ldots, s_n^{-1}r_n)$ for some elements $s_i\in S_{i, c}$ and $r_i\in R$. The element $r_i$ is unique up to adding an element of the ideal   $\ga_i$. Notice that $s_i\ga_i =0$ since $s_i\in S_{i,c}$, and so $s_i (r_i+\ga_i) = s_ir_i$. Then
\begin{eqnarray*}
 q&=&(s_1^{-2}s_1r_1, \ldots , s_n^{-2}s_nr_n) = \sum_{i=1}^n s_i^{-2}\cdot s_ir_i \\
 &=&(s_1^2+\cdots +s_n^2)^{-1}(s_1r_1+\cdots +s_nr_n) = c'^{-1} r
\end{eqnarray*}
where $c'=s_1^2+\cdots +s_n^2\in \CC'$ and $r= s_1r_1+\cdots +s_nr_n\in R$, since $S_{i,c}S_{j,c}=0$ for all $i\neq j$ (statement 6). Notice that $s_i^2\in S_{i,c}$ and $s_ir_i\in \bigcap_{j\neq i} \ga_j$ for $i=1, \ldots , n$ (by statement 6, since $s_i\in S_{i,c}$). Therefore,
$$Q'=\{\sum_{i=1}^n t_i^{-1}a_i\, | \, t_i\in S_{i, c}, a_i\in \bigcap_{j\neq i} \ga_j, i=1, \ldots , n\}.$$
This equality implies that $\CC'\in \Den_l(R, 0)$ and $\CC'^{-1}R= Q'$. Since $\CC'\subseteq \CC$ and $\CC \subseteq Q'^*$ (Step 1), every element $q\in Q'$ can be written as a left fraction $c'^{-1}r$ where $c'\in \CC'\subseteq \CC$ and $r\in R$. This fact implies that $ \CC \in \Ore_l(R)$: for given elements $c\in \CC$ and $r\in R$, $Q'\ni rc^{-1}=c'^{-1}r'$ for some elements $c'\in \CC'$ and $r'\in R$, or, equivalently, $c'r=r'c$, and so the left Ore condition holds for $\CC $.
  Since $\CC\in \Ore_l(R)$, $\CC \subseteq Q'^*$ and $R\subseteq Q'$, the map $\s' : Q=\CC^{-1}R\ra Q'$ is a ring monomorphism which is obviously an epimorphism as $Q' = \CC'^{-1}R$ and $\CC'\subseteq \CC$. Therefore, $\s'$ is an isomorphism and we can write  $Q= \CC^{-1}R=\CC'^{-1}R=Q'$, and statements 2 and 7 hold.

1. By Proposition \ref{A8Dec12}.(1) (or by Lemma \ref{1a27Nov12}),  $\CC \subseteq S_i$ for all $i=1, \ldots , n$, and so
\begin{eqnarray*}
 \CC &\subseteq & \bigcap_{i=1}^n S_i \stackrel{{\rm st.}\, 4}{=} \bigcap_{i=1}^nR\cap S_i'=R\cap \bigcap_{i=1}^n S_i'\\
 &=& R\cap R_1^*\times\cdots \times R_n^* = R\cap Q'^* = \CC,
\end{eqnarray*}
by Step 1. Therefore, $\CC  = \bigcap_{i=1}^nS_i$.

5. For all $i=1, \ldots , n$, $\emptyset \neq S_{i, c} \subseteq \bigcap_{j\neq i} \ga_j$ (statement 4), and so $\emptyset \neq S_{i,c} \subseteq S_i\backslash \bigcup_{j\neq i} S_j$, by statement 6. Let $r\in R$. Recall that $S_i= \s_i^{-1}(R_i^*)$, statement 4. Then $r\in S_i$  iff $\s_i(r)\in R_i^*$; and  $r\not\in S_i$ iff $\s_i(r) \in R_i\backslash R_i^* = R_i^0$. Therefore,
$$ S_i\backslash \bigcup_{j\neq i} S_j = R\cap R_i^* \times \prod_{j\neq i} R_j^0.$$

8. By statement 7, $Q= \{ \sum_{i=1}^n s_i^{-1}a_i \, | \, s_i\in S_{i,c}, \, a_i\in \bigcap_{j\neq i} \ga_j, \, i=1, \ldots , n\}$, and $Q=Q'$, by statement 2. Then $q\in Q^*=Q'^*$ iff $s_i^{-1}a_i\in R_i^*$ for $i=1, \ldots , n$ iff $a_i\in R_i^*$ for $i=1, \ldots , n$ iff $a_i\in S_{i,c}$ for $i=1, \ldots , n$ (since, by statement 6, $S_{i,c}= R_i^*\cap \bigcap_{j\neq i} \ga_j$). Therefore, $$q= (s_1+\cdots + s_n)^{-1} (a_1+\cdots +a_n)$$ where $s_1+\cdots + s_n, a_1+\cdots +a_n\in \CC'$, i.e. $Q^* = \{ s^{-1}t\, \ \, s,t\in \CC'\}$.  Since $\CC \subseteq Q^*$ (Step 1), the previous equality implies  the equality $Q^* = \{ s^{-1}t\, \ \, s,t\in \CC \}$.

9. By Step 1, $\CC = R\cap Q^*$. Then, by statement 8, $\CC =\{ s^{-1}t\, | \, s,t\in \CC'; s^{-1}t\in R\}$.
 $\Box $


\begin{corollary}\label{b3Jan13}
Suppose that the left quotient ring $Q$ of a ring $R$ is a simple left Artinian ring. Then $\maxDen_l(R) = \{ \CC \}$,  and so every left denominator set of the ring $R$ consists of regular elements.
\end{corollary}

{\it Proof}. Theorem \ref{31Dec12}.(1).  $\Box $


\begin{corollary}\label{a22Feb13}
 Let $R$ be a ring with a semisimple left quotient ring $Q$ (i.e. $R$ is a semiprime left Goldie ring) and we keep  the notation of Theorem \ref{31Dec12}. Then
\begin{enumerate}
\item The left localization radical of the ring $R$ is equal to zero and the set of regular elements of $R$ coincides with the set of completely left localizable elements.
\item
\begin{enumerate}
\item $\CN\CL_l(R) = \{ r\in R \, | \, \s_i (r)\in R_i^0  \; {\rm for}\; i=1, \ldots , n\}= \{ r\in R\, | \, r+\gp\not\in \CC_{R/\gp }\; {\rm for \; all}\;\; \gp \in \Min (R)\}$.
\item $R\cdot \CN\CL_l(R)\cdot R\subseteq \CN\CL_l(R)$.
\item $\CN\CL_l(R)+\CN\CL_l(R)\subseteq \CN\CL_l(R)$, i.e. $\CN\CL_l(R)$ is an ideal of  $R$, iff $\CN\CL_l(R) =0$ iff $R$ is a domain.
\end{enumerate}
\end{enumerate}
\end{corollary}

{\it Proof}. 1. Theorem \ref{31Dec12}.

2(a). The first equality follows from Theorem \ref{31Dec12}.(4). The second equality follows from Theorem \ref{1Jan13}(2c).

(b) The inclusion follows from the first equality in (a) and the fact that $R_iR_i^0R_i\subseteq R_i^0$ since $R_i$ is a simple Artinian ring.

(c) The statement (c) follows at once from the obvious fact that in a semisimple Artinian ring a sum of two zero divisors  is always a zero divisor iff the ring is a division ring. $\Box $

$\noindent $

The next corollary (together with Theorem \ref{31Dec12}) provides the necessary conditions of Theorem \ref{1Jan13}.

\begin{corollary}\label{b1Jan13}
 Let $R$ be as in Theorem \ref{31Dec12} and we keep the notation of Theorem \ref{31Dec12} and its proof. Then
\begin{enumerate}
\item The ring $R$ is a semiprime ring.
\item $\Min (R) = \{ \ga_1, \ldots , \ga_n\}$.
\item For each $i=1, \ldots , n$, $S_i = \pi_i^{-1} (\CC_{\bR_i}) = \{ c\in R\, | \, c+a_i\in \CC_{\bR_i}\}$ where $\pi_i : R\ra \bR_i : = R/ \ga_i$, $r\mapsto r+\ga_i$, and $\CC_{\bR_i}$ is the set of regular elements of the ring   $\bR_i$.
\item For each $i=1, \ldots , n$, $\CC_{\bR_i}=\overline{\s}_i^{-1}(R_i^*)\in \Ore_l(\bR_i)$ and $\CC_{\bR_i}^{-1} \bR_i\simeq R_i$  where $\overline{\s}_i:\bR_i\ra R_i$, $r+\ga_i\mapsto \frac{r}{1}$.
\item $S_i^{-1}\ga_j = \begin{cases}
0& \text{if }i=j,\\
R_i& \text{if }i\neq j.\\
\end{cases}$
\end{enumerate}
\end{corollary}

{\it Proof}. 1. Statement 1 follows from statement 2 and Theorem \ref{31Dec12}: $\bigcap_{\gp\in \Min (R)} \gp = \bigcap_{i=1}^n \ga_i =0$.

5. Clearly, $S_i^{-1}\ga_i =0$ for $i=1, \ldots , n$. Then statement 5 follows from Step 2 in the proof of Theorem \ref{31Dec12}.


2. Statement 2 follows from the following two statements:

(i) $\ga_1, \ldots , \ga_n$ {\em are prime ideals, }

(ii) $\Min (R) \subseteq \{ \ga_1, \ldots , \ga_n\}$.

Indeed, by (ii), $\Min (R) = \{ \ga_1, \ldots , \ga_m\}$, up to order.  The ideals $\ga_1, \ldots , \ga_n$ are incomparable, i.e. $\ga_i\not\subseteq \ga_j$ for all $i\neq j$. Then, by (i), we must have $m=n$ (since otherwise we would have (up to order)  $\ga_i\subseteq \ga_{m+1}$ for some $i$ such that $1\leq i\leq m$, a contradiction).

{\em Proof of (i)}: For each $i=1, \ldots , n$, the rings $R_i$ are left Artinian. So, if $\ga$ is an ideal of the ring $R_i$ then $S_i^{-1}\ga$ is an ideal of the ring $R_i$, hence $Q\ga = \prod_{i=1}^n S_i^{-1}\ga$ is an ideal of the ring $Q=Q'=\prod_{i=1}^nR_i$. If $\ga \gb \subseteq \ga_i$ for some ideals of the ring $R$ then $$Q\ga \cdot Q\gb\subseteq Q\ga \gb \subseteq Q\ga_i = \prod_{j=1}^n S_j^{-1}\ga_i = R_1\times \cdots \times 0 \times \cdots \times R_n = \ga_i',$$
by statement 5 and Theorem \ref{31Dec12}.(3). The ideal $\ga_i'$ of the ring $Q$ is a prime ideal. Therefore, either $Q\ga \subseteq \ga_i'$ or $Q\gb \subseteq \ga_i'$. Then, either $\ga \subseteq R\cap Q\ga \subseteq R\cap \ga_i'= \ga_i$ or $\gb \subseteq R\cap Q\gb \subseteq R\cap \ga_i'= \ga_i$ (since $\ga_i = R\cap \ga_i'$, by Theorem \ref{31Dec12}.(4)), i.e. $\ga_i$ is a prime ideal of the ring $R$.

{\em Proof of (ii)}: Let $\gp$ be a prime ideal of the ring $R$. Then $$\prod_{i=1}^n \ga_i \subseteq \bigcap_{i=1}^n \ga_i =\ker (\s )= \{ 0 \}\subseteq \gp$$ (Theorem \ref{31Dec12}), and so $\ga_i \subseteq \gp$ for some $i$, and the statement (ii) follows from the statement (i).

3 {\rm and} 4. There is a commutative diagram of ring homomorphisms

$$\xymatrix{%
R\ar[rd]^{\pi_i}\ar[r]^{\s_i} & R_i\\
 &\bR_i\ar[u]_{\overline{\s}_i}}
$$
where $\overline{\s}_i: r+\ga_i\mapsto\frac{r}{1}$ is a monomorphism, $\pi_i(S_i) \in\Den_l(\bR_i, 0)$ and $ \pi_i(S_i)^{-1}\bR_i \simeq R_i=S_i^{-1}R$ (via an obvious extension of $\overline{\s}_i$). The ring $R_i$ is a simple left Artinian ring, hence a left localization maximal ring. By Theorem \ref{15Nov10}.(3), $S_i=\s_i^{-1}(R_i^*)$. Clearly, $\overline{\s}_i^{-1}(R_i^*)\subseteq \CC_{\bR_i}$. On the other hand, for each element $c\in \CC_{\bR_i}$, the map $\cdot c: R_i\ra R_i$, $r\mapsto rc$, is a left $R_i$-module monomorphism, hence as isomorphism, and so $ c\in R_i^*$, i.e. $\CC_{\bR_i}\subseteq R_i^*$. Therefore, $\overline{\s}_i^{-1}(R_i^*)= \CC_{\bR_i}$. Now,
\begin{equation}\label{SsiR}
 S_i=\s^{-1}(R_i^*) = (\overline{\s}_i\pi_i)^{-1}(R_i^*)=
 \pi_i^{-1}(\overline{\s}_i^{-1}(R_i^*))=\pi_i^{-1}(\CC_{\bR_i}).\;\;\;
\end{equation}
Therefore, $\CC_{\bR_i} =\pi_i(S_i)\in \Ore_l(\bR_i)$ and $\CC_{\bR_i}^{-1} \bR_i \simeq \pi_i(S_i)^{-1}\bR_i\simeq R_i$. $\Box$



\section{The Second  Criterion   (via the minimal primes)}\label{TSECMP}

The aim of this section is to give another criterion (Theorem \ref{1Jan13}) for a ring $R$ to have a semisimple left quotient ring $Q$.

\begin{theorem}\label{1Jan13}
 ({\bf The Second Criterion})  Let $R$ be a ring. The following statements are equivalent.
\begin{enumerate}
\item The ring $R$ has a semisimple left quotient ring $Q$.
\item
\begin{enumerate}
\item The ring $R$ is a semiprime ring.
\item The set $\Min (R) $ of minimal primes of the ring $R$ is a finite set, say,  $\{ \ga_1, \ldots , \ga_n\}$.
    \item For each $i=1, \ldots , n$, $S_i := \pi_i^{-1}(\CC_{\bR_i}) = \{ c\in R\, | \, c+\ga_i \in \CC_{\bR_i}\}\in \Den_l(R, \ga_i)$ where $\pi_i : R\ra \bR_i:=R/ \ga_i$, $r\mapsto r+\ga_i$.
        \item For each $i=1, \ldots , n$, the ring $R_i:=S_i^{-1}R$ is a simple left Artinian ring.
\end{enumerate}
\end{enumerate}
If one of the two equivalent conditions holds then $\maxDen_l(R)=\{ S_1, \ldots , S_n\}$ and $Q\simeq \prod_{i=1}^n R_i$.
\end{theorem}

{\it Proof}. $(1\Rightarrow 2)$ Corollary \ref{b1Jan13} and Theorem \ref{31Dec12}.

$(2\Rightarrow 1)$ It suffices to prove the following claim.

{\em Claim}: $\maxDen_l(R)=\{ S_1, \ldots , S_n\}$.

Since then $\bigcap_{i=1}^n \ga_i =0$ (as $R$ is a semiprime ring) and so the assumptions of Theorem \ref{31Dec12} hold and as a result the ring $R$ has a semisimple left quotient ring.

{\it Proof of the Claim}:

(i) $S_1, \ldots , S_n$ {\em are distinct} since $\ga_1= \ass (S_1), \ldots , \ga_n= \ass (S_n)$ are distinct.

(ii) $\maxDen_l(R)\supseteq \{ S_1, \ldots , S_n\}$: Every simple left Artinian ring is a left localization maximal ring, for example, $R_i= S_i^{-1}R\simeq \pi_i(S_i)^{-1}\bR_i$ for $i=1, \ldots , n$. The ring $R_i$ is a left Artinian ring that contains the ring $\bR_i$ (via $\overline{\s}_i$). Hence, $\CC_{\bR_i}\subseteq R_i^*$. Since $\pi_i (S_i) \subseteq \CC_{\bR_i}$, $\pi_i (S_i) \in \Ore_l(\bR_i)$, $\CC_{\bR_i} \subseteq R_i^*$ and $R_i =\pi_i(S_i)^{-1} \bR_i$, we see that $$\CC_{\bR_i} \in \Ore_l(\bR_i)\;\; {\rm  and }\;\; \CC_{\bR_i}^{-1} \bR_i = \pi_i(S_i)^{-1} \bR_i = R_i.$$
 Recall that $\s_i : R\ra R_i$, $r\mapsto \frac{r}{1}$, and $\pi_i: R\ra \bR_i$, $r\mapsto \br = r+\ga_i$. Clearly, $\s_i = \overline{\s}_i\pi_i$. The ring $\CC_{\bR_i}^{-1}\bR_i \simeq R_i$ is a left localization maximal ring. Therefore, $\maxDen_l(\bR_i)=\{ \CC_{\bR_i}\}$, by Proposition \ref{A8Dec12}, and $\CC_{\bR_i}= \overline{\s}_i^{-1}(R_i^*)$, by Theorem \ref{15Nov10}.(3). Then, by Theorem \ref{15Nov10}.(3), $\s_i^{-1} (R_i^*)\in \maxDen_l(R_i)$ since $R_i$ is a left localization maximal ring. Now,
$$\s_i^{-1} (R_i^*) = (\overline{\s}_i\pi_i)^{-1}(R_i^*) = \pi_i^{-1}(\overline{\s}_i^{-1}(R_i^*))=\pi_i^{-1} (\CC_{\bR_i})= S_i.$$
Then $S_i\in \maxDen_l(R$ for $i=1, \ldots , n$.

(iii) $\maxDen_l(R)\subseteq \{ S_1, \ldots , S_n\}$: We have to show that, for a given left denominator set $S\in \Den_l(R, \ga )$ of the ring $R$; $S\subseteq S_i$ for some $i$. We claim that $S\cap \ga_i =\emptyset$ for some $i$ otherwise for each $i=1, \ldots , n$, we would have chosen an element $s_i\in S\cap \ga_i$, then we would have $$\prod_{i=1}^n s_i\in S\cap \bigcap_{i=1}^n \ga_i = S\cap 0=\emptyset ,$$ a contradiction. We aim to show that $S\subseteq S_i$. The condition $S\cap \ga_i = \emptyset$ implies that $\bS :=\pi_i (S) \in \Ore_l(\bR_i)$. Let $\gb = \ass_{\bR_i}(S)$.

($\alpha$) We claim that $\gb =0$. The ring $R_i$ is a simple left Artinian ring with $\bR_i\subseteq R_i$. In particular, it satisfies  the a.c.c. on right annihilators  (since $R_i$ is a matrix ring with entries from a division ring). By Lemma \ref{a1Jan13} which is applied in the situation that $\bR_i\subseteq R_i$ and $\bS\in \Ore_l(\bR_i)$, the core $\bS_c$ of the left Ore set $\bS$ is a non-empty set. Let $s\in \bS_c$. Then $\bR_i s\bR_i\cdot \gb = \bR_i s\gb =0$. The ring $\bR_i$ is a prime ring and $\bR_i s\bR_i\neq 0$ since $0\neq s\in \bR_i s\bR_i$. Therefore, $\gb =0$.

($\beta $) $\ga \subseteq \ga_i$: by ($\alpha $).

($\g $) $SS_i\in \Den_l(R)$, by ($\beta $) and Lemma \ref{1a27Nov12}.    Hence $S\subseteq SS_i= S_i$, by the maximality  of $S_i$.  $\Box $


The next result is a useful tool in finding elements of the core of an Ore set in applications.

\begin{lemma}\label{a1Jan13}
Let $R\subseteq R'$ be rings and $S\in \Ore_l(R)$. Suppose that the ring $R'$ satisfies the a.c.c. on right annihilators and let $\Max_{R'}(S)=\{ s\in S\, | \, \ker_{R'}(s\cdot ) $ is a maximal element of the set $\{ \ker_{R'} (s\cdot)\, \ \, s\in S\}\, \}$. Then
\begin{enumerate}
\item $\Max_{R'}(S) \neq \emptyset$ and $\ker_{R'}(s\cdot ) = \ker_{R'}(t\cdot )$ for all elements $s,t\in \Max_{R'}(S)$.
\item $\Max_{R'}(S)\subseteq S_c$.
\end{enumerate}
\end{lemma}

{\it Proof}. 1. The ring $R'$ satisfies the a.c.c. on right annihilators, so  $\Max_{R'}(S) \neq \emptyset$. Since $S\in \Ore_l(R)$, we must have $\ker_{R'}(s\cdot ) = \ker_{R'}(t\cdot )$ for all elements $s,t\in \Max_{R'}(S)$ (the left Ore condition for $S$ implies that $s's=rt$ for some elements $s'\in S$ and $r\in R$, hence $s's\in S$ and $\ker_{R'}(s\cdot ), \ker_{R'}(t\cdot )\subseteq \ker_{R'}(ss'\cdot )$, and so
$$\ker_{R'}(s\cdot )=\ker_{R'}(ss'\cdot )= \ker_{R'}(t\cdot ) \,). $$

2. For all $s\in S$, $\ker_{R}(s\cdot )=R\cap \ker_{R'}(s\cdot )$. Therefore, statement 2 follows from statement 1.   $\Box $



\section{The Third Criterion (in the spirit of Goldie-Lesieur-Croisot)
}\label{THIGLCR}

The aim of this section is to prove Theorem \ref{20Jan13}. This is the Third Criterion for a ring to have a semisimple left quotient ring. It is close to to Goldie's Criterion but in applications it is easier to check its conditions. It reveals the `local' nature of the fact that a ring has a semisimple left quotient ring.

 For a semiprime  ring $R$ and its ideal $I$, the left annihilator of $I$ in $R$ is equal to the right annihilator of $I$ in $R$ and is denoted $\ann (I)$. A ring is called a {\em left Goldie ring} if it has a.c.c. on left annihilators  and does not contain infinite direct sums of nonzero left ideals.

\begin{theorem}\label{20Jan13}
({\bf The Third Criterion}) Let $R$ be a ring. The following statements are equivalent.
\begin{enumerate}
\item The ring $R$ has a semisimple left quotient ring.
\item The ring $R$ is a semiprime ring with $|\Min (R)|<\infty$ and, for each $\gp \in \Min (R)$, the ring $R/\gp$ is a left Goldie ring.
\end{enumerate}
\end{theorem}

{\it Remark}. The condition $|\Min (R)|<\infty$ in Theorem \ref{20Jan13} can be replaced by any of the  equivalent conditions of Theorem \ref{2.2.15-[MR]}.

\begin{theorem}\label{2.2.15-[MR]}
{\rm (Theorem 2.2.15, \cite{MR})} The following conditions on a semiprime ring $R$ are equivalent.
\begin{enumerate}
\item ${}_RR_R$ has finite uniform dimension.
\item $|\Min (R)|<\infty$.
\item $R$ has finitely many annihilator ideals.
\item $R$ has a.c.c. on annihilator ideals.
\end{enumerate}
\end{theorem}

{\bf Proof of Theorem \ref{20Jan13}}. $(1\Rightarrow 2)$ Theorem \ref{1Jan13} and Goldie's Theorem (since $R_i =S_i^{-1}R\simeq \CC_{R/\ga_i}^{-1}(R/\ga_i) = Q(R/\ga_i )$).

$(1\Leftarrow 2)$  If the ring $R$ is a prime ring then the result follows from Goldie-Lesieur-Croisot's Theorem  in the prime case \cite{Goldie-PLMS-1958}, \cite{Lesieur-Croisot-1959}. So, we can assume that $R$ is  not a prime ring, that is $n:=|\Min (R)|\geq 2$.
 The idea of the proof of the implication  is to show that the conditions (a)-(d) of Theorem \ref{1Jan13} hold. The conditions (a) and (b) are given.  Let $\min (R) = \{ \ga_1, \ldots , \ga_n\}$, we keep the notation of Theorem \ref{1Jan13}. For each $i=1, \ldots , n$, the ring $\bR_i = R/ \ga_i$ is a prime left Goldie ring (by the assumption). By Goldie-Lesieur-Croisot's Theorem  in the prime case,  $\CC_{\bR_i}\in \Ore_l(\bR_i)$ and  $Q(\bR_i)$ is a simple Artinian ring. The conditions (c) and (d) follows from Proposition \ref{a20Jan13} and the Claim.

{\it Claim}: $\ann (\ga_i) \cap S_i\neq \emptyset$ for $i=1, \ldots ,n$.

Indeed, by Proposition \ref{a20Jan13} which is applied in the situation  $R$, $I=\ga_i$, $\overline{\ga}=0$ and $\bS = \CC_{\bR_i}$, we have $S_i\in \Den_l(R, \ga_i )$ and $S_i^{-1}R \simeq \CC_{\bR_i}^{-1}\bR_i = Q(\bR_i)$ is a simple Artinian ring, as required.

{\it Proof of the Claim}. The ring $R$ is  semiprime, so $\ann (\ga_i) +\ga_i = \ann (\ga_i) \oplus \ga_i$ (since $(\ann (\ga_i) \cap \ga_i)^2=0$, we must have $\ann (\ga_i) \cap \ga_i =0$, as $R$ is a semiprime ring).
 Since $n\geq 2$, $0\neq \bigcap_{j\neq i}\ga_i\subseteq \ann (\ga_i)$ (as $\ann (\ga_i) \cdot \bigcap_{j\neq i}\ga_i\subseteq \bigcap_{k=1}^n\ga_k=0$, the ring $R$ is a semiprime ring). Let $\pi_i: R\ra \bR_i:= R/\ga_i$, $r\mapsto r+\ga_i$.
Then $\pi_i(\ann (\ga_i))$ is a {\em nonzero} ideal of the prime ring $\bR_i$. Any nonzero ideal of a prime ring is an essential left (and right) ideal, Lemma 2.2.1.(i), \cite{MR}.  Hence, $\pi (\ann (\ga_i))\cap \CC_{\bR_i}\neq \emptyset$, by Goldie-Lesieur-Croisot's Theorem (in the prime case). Therefore, $\ann (\ga_i) \cap S_i\neq \emptyset$, as required. $\Box $


\begin{proposition}\label{a20Jan13}
Let $R$ be a ring, $I$ be its ideal, $\pi : R\ra \bR :=R/I$, $r\mapsto \br =r+I$, $\bS\in \Den_l(\bR , \overline{\ga}  )$, $S:=\pi^{-1} (\bS )$ and $\ga :=\pi^{-1}(\overline{\ga}  )$. If, for each element $x\in I$, there is an element $s\in S$ such that $sx=0$ then $S\in \Den_l(R, \ga )$ and $ S^{-1}R\simeq \bS^{-1}\bR$.
\end{proposition}

{\it Proof}. By the very definition, $S$ is a multiplicative set.

(i) $S\in \Ore_l(R)$: Given elements $s\in S$ and $r\in R$, we have to find elements $s'\in S$ and $r'\in R$ such that $s'r=r's$. Since $\bS\in \Den_l(\bR , \overline{\ga} )$, $\bs_1\br = \br_1\bs $ for some elements $ s_1\in S$ and $r_1\in R$. Then $x:= s_1r-r_1s\in I$, and so $0=s_2x$ for some element $s_2\in S$, then  $s_2s_1r=s_2r_1s$. It suffices to take $s'=s_2s_1$ and $ r_1'= s_2r_1$.

(ii) $\ass (S) = \ga$: Let $r\in \ass (S)$, i.e. $sr=0$ for some element $s\in S$. Then $\bs \br =0$, and so $\br \in \overline{\ga} $. This implies that $r\in \ga$. So, $\ass (S) \subseteq \ga$. Conversely, suppose that $a\in \ga$. We have to find an element $s\in S$ such that $sa=0$. Since $\overline{a} \in \overline{\ga} $ and $\bS\in \Den_l(\bR , \overline{\ga} )$, $\bs_1 \overline{a}=0$ for some element $s_1\in S$. Then $s_1a\in I$, and so $s_2s_1a=0$ for some element $s_2\in S$. It suffice to take $s=s_2s_1$.

(iii) $S\in \Den_l(R, \ga )$: In view of (i) and (ii), we have to show that if $rs=0$ for some $r\in R$ and $s\in S$ then $r\in \ga$. The equality  $\br \bs =0$ implies that $\br \in \overline{\ga}$   (since $\bS\in \Den_l(\bR , \overline{\ga} )$), and so $r\in \ga$.

(iv) $S^{-1}R\simeq \bS^{-1}\bR$: By the universal property of left localization, the map $S^{-1}R\ra \bS^{-1}\bR$, $s^{-1}r\mapsto \bs^{-1}\br$, is a well-defined ring homomorphism which is obviously an epimorphism.  It suffices to show that the kernel is zero. Given $s^{-1}r\in S^{-1}R$ with $\bs^{-1}\br =0$. Then $\bs_1\br =0$ in $\bR$, for some element $s_1\in S$. Then $s_1r\in I$, and so $s_2s_1r=0$ for some element $s_2\in S$. This means that $ s^{-1}r=0$ in $S^{-1}R$. So, the kernel is equal to zero and as the result   $S^{-1}R\simeq \bS^{-1}\bR$. $\Box $


\section{The Fourth Criterion (via certain left denominator sets)
}\label{TFCCLD}

 In this section, a very useful criterion for  a ring $R$ to have a semisimple left quotient ring is given (Theorem \ref{22Jan13}). It is a corollary of Theorem \ref{20Jan13} and the following characterization of the set of minimal primes in  a semiprime ring.

\begin{proposition}\label{A22Jan13}
Let $R$ be a ring and $\ga_1, \ldots , \ga_n$ ideals of $R$ and $n\geq 2$. The following statements are equivalent.
\begin{enumerate}
\item The ring $R$ is a semiprime ring with $\Min (R) = \{ \ga_1, \ldots , \ga_n\}$.
\item The ideals $\ga_1, \ldots , \ga_n$ are incomparable prime ideals  with $\bigcap_{i=1}^n \ga_i=0$.
\end{enumerate}
\end{proposition}

{\it Proof}. $(1\Rightarrow 2)$ Trivial.

$(1\Leftarrow 2)$  The ideals  $\ga_1, \ldots , \ga_n$ are incomparable ideals hence are distinct. The ideals $\ga_i $ are  prime.  So, or each $i=1, \ldots , n$, fix a minimal prime ideal $\gp_i\in \Min (R)$  such that $\gp_i\subseteq \ga_i$. We claim that $\gp_i=\ga_i$ for $i=1, \ldots , n$. The inclusions $\prod_{j=1}^n\ga_j\subseteq \bigcap_{j=1}^n \ga_j=\{0\} \subseteq \gp_i$, imply $\ga_j\subseteq \gp_i\subseteq \ga_i$ for some $j$, necessarily $j=i$ (since the ideals $\ga_1, \ldots , \ga_n$ are incomparable). Therefore, $\gp_i=\ga_i$. So, the ideals $\ga_1, \ldots , \ga_n$ are minimal primes of $R$ with $\bigcap_{i=1}^n \ga_i=0$. Therefore, $\Min (R) = \{ \ga_1, \ldots , \ga_n\}$ (otherwise we would have a minimal prime $\gp\in \Min (R)$ distinct from the ideals $\ga_1, \ldots , \ga_n$ but then $\prod_{i=1}^n \ga_i\subseteq \bigcap_{i=1}^n \ga_i = \{ 0 \} \subseteq \gp \Rightarrow \ga_i\subseteq \gp$, for some $i$,  a contradiction).  $\Box $


\begin{theorem}\label{22Jan13}
({\bf The Fourth Criterion}) Let $R$ be a ring. The following statements are equivalent.
\begin{enumerate}
\item The ring $R$ has a semisimple left quotient ring.
\item There are left denominator sets $S_1', \ldots , S_n'$ of the ring $R$ such that the rings $R_i:=S_i'^{-1}R$, $i=1, \ldots , n$, are simple left Artinian rings and the map
    $$ \s := \prod_{i=1}^n \s_i : R\ra \prod_{i=1}^n R_i, \;\; R\mapsto (\frac{r}{1}, \ldots , \frac{r}{1}), $$
    is an injection where $\s_i : R\ra R_i$, $r\mapsto \frac{r}{1}$.
\end{enumerate}
If one of the equivalent conditions holds then the set $\maxDen_l(R)$ contains precisely the distinct elements of the set $\{  \s_i^{-1}(R_i^*)\, | \, i=1, \ldots , n\}$.
\end{theorem}

{\it Proof}. $(1\Rightarrow 2)$ By Theorem \ref{31Dec12} and (\ref{llRseq}), it suffices to take $\maxDen_l(R)$ since,  by (\ref{llradR}), $\gll_R=0$.

$(1\Leftarrow 2)$ If $n=1$ the implication $(2\Rightarrow 1)$ is obvious: the map $\s : R\ra S_1'^{-1}R$ is an injection, hence $S_1'\in \Den_l(R, 0)$. In particular, the set $S_1$ consists of {\em regular} elements of the ring $R$. Since $R\subseteq S_1'^{-1}R$ (vis $\s$) and the ring $S_1'^{-1}R$ is a simple left  Artinian ring, $S_1'\subseteq \CC\subseteq (S_1'^{-1}R)^*$, and so  $\CC \in \Ore_l(R)$ and $Q=S_1^{-1}R$.

Suppose that $n\geq 2$. It suffices to verify that the conditions of Theorem \ref{20Jan13}.(2) hold. For each $i=1, \ldots , n$, the ring $R_i$ is a simple left Artinian ring, and so $R_i$ is a left localization maximal ring. By Theorem \ref{15Nov10}.(3),
 $$S_i=\s_i^{-1}(R_i^*)\in \maxDen_l(R), \;\;\; S_i^{-1}R\simeq R_i,$$
and obviously $S_i'\subseteq S_i$ and $\ga_i= \ass (S_i ) =\ass (S_i')$. Up to order, let $S_1, \ldots , S_m$ be the distinct elements of the collection $\{ S_1, \ldots , S_n\}$. By Proposition \ref{b27Nov12}, the ideals $\ga_1, \ldots , \ga_m$ are incomparable. The map $\s$ is an injection, hence
$$ 0=\ker (\s ) = \bigcap_{i=1}^n \ass (S_i')=\bigcap_{i=1}^n \ass (S_i)=\bigcap_{i=1}^m \ass (S_i').$$
 For each $i=1, \ldots , m$, let $\pi_i : R\ra \bR_i:=R/ \ga_i $, $r\mapsto \br_i = r+\ga_i$. Then $\bS_i:=\pi_i(S_i)\in \Den_l(\bR_i, 0)$ and the ring $\bS_i^{-1}\bR_i\simeq S_i^{-1}R$ is a simple left Artinian ring which necessarily coincides with the left quotient ring $Q(\bR_i)$ of the ring $\bR_i$ (see the case $n=1$). By Goldie-Lesieur-Croisot's Theorem the ring $\bR_i$ is a prime left Goldie ring. By Proposition \ref{A22Jan13}, $R$ is a semiprime ring with $\Min (R) = \{ \ga_1, \ldots , \ga_m \}$. So, the ring $R$ satisfies the conditions of Theorem \ref{20Jan13}.(2), as required. $\Box $



\section{Left denominator sets of finite direct products of rings}\label{LDSDFPR}

Let $R=\prod_{i=1}^nR_i$ be a direct product of rings $R_i$. The aim of this section is to give a criterion  of when a multiplicative subset of $R$ is a left Ore/denominator set (Proposition \ref{b27Dec12}) and using it to show that the set of maximal left denominator sets of the ring $R$ is the union of the sets of maximal left denominators sets of the rings $R_i$ (Theorem \ref{c26Dec12}).

Let $R=\prod_{i=1}^nR_i$ be a direct product of rings $R_i$ and $1=\sum_{i=1}^ne_i$ be the corresponding sum of central orthogonal idempotents ($e_i^2=e_i$ and $e_ie_j=0$ for all $i\ne j$). Let $\pi_i: R\ra R_i$, $r=(r_1, \ldots , r_n)\mapsto r_i$ and $\nu_i:R_i\ra R$, $r_i\mapsto (0, \ldots , 0,r_i, 0, \ldots , 0)$. Clearly, $\nu_i\pi_i=\id_{R_i}$, the identity map in $R_i$. Every ideal $\ga$ of the the ring $R$ is the direct product  $\ga= \prod_{i=1}^n \ga_i$ of ideals $\ga_i= e_i\ga$ of the rings $R_i$.

\begin{lemma}\label{a27Dec12}
Let $R=\prod_{i=1}^iR_n$ be a direct product of rings $R_i$. We keep the notation as above.
\begin{enumerate}
\item If $S\in \Ore_l(R, \ga = \prod_{i=1}^n\ga_i)$ then for each $i=1, \ldots , n$ either $0\in S_i:=\pi_i(S)$ or otherwise $S_i\in \Ore_l(R_i, \ga_i)$ and at least for one $i$,  $S_i\in \Ore_l(R_i, \ga_i)$.
\item If $S\in \Den_l(R, \ga = \prod_{i=1}^n\ga_i)$ then for each $i=1, \ldots , n$ either $0\in S_i:=\pi_i(S)$ or otherwise $S_i\in \Den_l(R_i, \ga_i)$ and at least for one $i$, $S_i\in \Den_l(R_i, \ga_i)$.
\end{enumerate}
\end{lemma}

{\it Proof}.  Straightforward. $\Box $

$\noindent $

For each $S\in \Ore_l(R)$, let $\supp (S):=\{ i\, | \, \pi_i(S)\in \Ore_l(R_i)\}$. By Lemma \ref{a27Dec12}.(1),
\begin{equation}\label{aprOre}
\ga = \ass (S) = \prod_{i=1}^n\ga_i, \;\; \ga_i = \begin{cases}
\ass_{R_i}(\pi_i(S))& \text{if }i\in \supp (S),\\
R_i& \text{if }i\not\in \supp (S).\\
\end{cases}
\end{equation}
The following  proposition  is a criterion for a multiplicative set of a finite direct product of rings to be a left Ore set or a left denominators set.

\begin{proposition}\label{b27Dec12}
Let $R=\prod_{i=1}^nR_i$ be a direct product of rings $R_i$, $S$ be a multiplicative set of $R$ and $S_i:=\pi_i(S)$ for $i=1, \ldots , n$.  Then
\begin{enumerate}
\item $S\in \Ore_l(R)$ iff for each $i=1, \ldots , n$
either $0\in S_i$ or otherwise $S_i\in \Ore_l(R_i)$ and at least for one $i$,  $S_i\in \Ore_l(R_i)$.
\item $S\in \Den_l(R)$ iff for each $i=1, \ldots , n$ either $0\in S_i$ or otherwise $S_i\in \Den_l(R_i)$ and at least for one $i$, $S_i\in \Den_l(R_i)$.
\end{enumerate}
\end{proposition}

{\it Proof}. 1. $(\Rightarrow )$ Lemma \ref{a27Dec12}.(1).

$(\Leftarrow )$ Clearly, $\supp (S)\neq \emptyset$ and $\ass (S)$ is given by (\ref{aprOre}).
 Without loss of generality we may assume that $\supp (S) = \{ 1, \ldots , m\}$ where $1\leq m \leq n$. We have to show that for given elements $s=(s_i)\in S$ and $r=(r_i)\in R$ there are elements $s'=(s'_i)\in S$ and $r'=(r'_i)\in R$ such that $s'r=r's$. Since $S_1\in \Ore_l(R_1)$, $t_{11} r_1=a_1s_1$ for some elements $t_{11}\in S_1$ and $a_1\in R_1\subseteq R$. Fix an element $t_1\in S$ such that $\pi_1(t_1) = t_{11}$. Using the same sort of argument but for $S_2\in \Ore_l(R_2)$, we have $ t_{22}\pi_2(t_1r) = a_2s_2$ for some elements $t_{22}\in S_2$ and $a_2\in R_2\subseteq R$. Fix an element $t_2\in R$ such that $\pi_2(t_2)=t_{22}$. Repeating the same trick several times  we have the equalities
 $$ t_{ii}\pi_i(t_{i-1}\cdots t_2t_1r) = a_is_i, \;\; i=1, \ldots , m, $$ where $t_{ii}\in S_i$, $a_i\in R_i$ and $t_i\in S$ such that $\pi_i(t_i) = t_{ii}$.

 If $m=n$ then it suffices to take $s'=t_n\cdots t_2t_1$ and the element $r'=(r_i')$ can be easily written using the equalities above
  $$ r_1'=\pi_1(t_nt_{n-1}\cdots t_2)a_1, \; r_2'=\pi_2(t_nt_{n-1}\cdots t_3)a_2,\; \ldots , \; r_{n-1}'=\pi_{n-1}(t_n)a_{n-1}, \; r_n'=a_n.$$
If $m<n$ then this case can deduced to the previous one. For, each $j=m+1, \ldots , n$, choose an element $\th_j=(\th_{ji})\in S$ with $\th_{jj}=0$. Let $\th = \th_{m+1}\cdots \th_n$. It suffices to take $s'=\th t_mt_{m-1}\cdots t_1$ and $r'= \th \cdot (r_1'', \ldots , r_m'', 0, \ldots , 0)$ where
$$ r_1''=\pi_1(t_mt_{m-1}\cdots t_2)a_1, \; r_2''=\pi_2(t_mt_{m-1}\cdots t_3)a_2,\; \ldots , \; r_{m-1}''=\pi_{m-1}(t_m)a_{m-1}, \; r_m''=a_m.$$

2. $(\Rightarrow )$  Lemma \ref{a27Dec12}.(2).

$(\Leftarrow )$ By statement 1, $S\in \Ore_l(R)$. Without loss of generality we may assume that $\supp (S) = \{ 1, \ldots , m\}$ where $1\leq m \leq n$. By (\ref{aprOre}), $\ga= \prod_{i=1}^m\ga_i\times \prod_{j=m+1}^nR_j$. Using the fact that $S_i\in \Den_l(R_i, \ga_i)$ for $i=1, \ldots , m$ and $\ker (\th\cdot ) \supseteq \prod_{j=m+1}^nR_j$, we see that $S\in \Den_l(R, \ga )$: If $rs=0$ for some elements $r=(r_i)\in R$ and $s= (s_i)\in S$ then $r_i\in \ga_i$ for $i=1, \ldots , m$ since $S_i\in \Den_l(R_i, \ga_i)$. Fix an element $t_1\in S$ such that $\pi_1(t_1r)=0$. Since $\pi_2(t_1r)\in \ga_2$, fix an element $t_2\in S$ such that $\pi_2(t_2t_1r)=0$. Repeating the same argument several times we find elements $t_1, \ldots , t_m\in S$ such that
$$ \pi_i(t_it_{i-1}\cdots t_1r)=0\;\; {\rm for }\;\; i=1, \ldots , m.$$
Then $s'r=0$ where $s'=\th t_it_{i-1}\cdots t_1\in S$ and the element $\th$ is defined in the proof of statement 1.  $\Box $

$\noindent $

By Proposition \ref{b27Dec12}, for each $i=1, \ldots , n$, we have the injection
\begin{equation}\label{aab1}
\maxDen_l(R_i) \ra \maxDen_l(R=\prod_{i=1}^nR_i), \;\; S_i\mapsto R_1\times \cdots \times S_i\times \cdots \times R_n.
\end{equation}
We identify $\maxDen_l(R_i)$ with its image in  $\maxDen_l(R=\prod_{i=1}^nR_i)$.


$\noindent $

{\bf  Proof of Theorem \ref{c26Dec12}}. The theorem follows at once from Proposition \ref{b27Dec12}.(2). $\Box $

$\noindent $

\begin{corollary}\label{bc27Dec12}
Let $R$ be a  ring and $1=e_1+\cdots +e_n$ be a sum of central orthogonal idempotents of the ring $R$. Then for each  $S\in \maxDen_l(R)$ precisely one idempotent belongs to $S$.
\end{corollary}

{\it Proof}. The ring $R=\prod_{i=1}^nR_i$ is a direct product of rings $R_i:=e_iR$. Then the result follows from Theorem \ref{c26Dec12}.  $\Box $



\section{Criterion for $R/\gll_R$ to have a semisimple left quotient  ring}\label{CRLSLG}

The aim of this section  is to give  a criterion for the factor ring $R/ \gll$ (where $\gll$ is the left localization radical of $R$) to have a semisimple left quotient ring (Theorem \ref{24Feb13}), i.e. $R/ \gll$ is a semiprime left Goldie ring. Its proof is based on four criteria and some results (Lemma \ref{a24Feb13} and Lemma \ref{b24Feb13}).

\begin{theorem}\label{24Feb13}
Let $R$ be a ring, $\gll = \gll_R$ be the left localization radical of $R$, $\Min (R, \gll )$ be the set of minimal primes over the ideal $\gll$ and $\pil : R\ra R/ \gll$, $ r\mapsto \br =r+\gll$. The following statements are equivalent.
\begin{enumerate}
\item The ring $R/ \gll$ has a semisimple left quotient ring, i.e. $R/ \gll$ is a semiprime left Goldie ring.
\item
\begin{enumerate}
\item $|\maxDen_l(R)|<\infty$.
\item For every $S\in \maxDen_l(R)$, $S^{-1}R$ is a simple left Artinian ring.
\end{enumerate}
\item
\begin{enumerate}
\item $\gll = \bigcap_{\gp\in \Min (R, \gll )}\gp$.
\item $\Min (R, \gll )$ is a finite set.
\item For all $\gp\in  \Min (R, \gll )$, the set $S_\gp := \{ c\in R\, | \, c+\gp \in \CCRp \}$ is a left denominator set of the ring $R$ with $\ass (S_\gp ) = \gp$.
    \item For each $\gp\in \Min (R, \gll )$, $S_\gp^{-1}R$ is a simple left Artinian ring.
        \item For each $\gp\in \Min (R, \gll )$ and for each $l\in \gll$, there is an element $s\in S_\gp$ such that $sl=0$.
\end{enumerate}
\item
 \begin{enumerate}
\item $\gll = \cap_{\gp \in \Min (R, \gll )}\gp$.
\item $\Min (R, \gll )$ is a finite set.
\item For each $\gp\in \Min (R, \gll )$, $R/ \gp$ is a left Goldie ring.
\end{enumerate}
\end{enumerate}
If one of the equivalent conditions 1--3 holds then

(i) The map $\pil': \maxDen_l(R)\ra \maxDen_l(R/\gll )$, $ S\mapsto \pil (S)$, is a bijection with the inverse $T\mapsto \pil^{-1}(T)$.

(ii) $\maxDen_l(R)= \{ S_\gp \, | \, \gp \in \Min (R, \gll )\}$.

(iii) $\maxDen_l(R/ \gll )\{ S_{\gp / \gll}\, | \, \gp \in \Min (R, \gp )\}$ where $S_{\gp / \gll}:=\pil (S_\gp )$.

(iv) For all $\gp \in \Min (R)$,  $S_\gp = \pil^{-1}(S_{\gp / \gll})$, $\ass (S_{\gp / \gll}) = \ass (S_\gp )/ \gll$ and $ S^{-1}_\gp R\simeq S_{\gp / \gll}^{-1} R/ \gll$ is a simple left Artinian ring.
\end{theorem}

Theorem \ref{24Feb13} shows that when the factor ring $R/ \gll$ is a semiprime left Goldie ring there are tight connections between the localization properties of the rings $R$ and $R/ \gll$. In particular, the map $\pil'$ is a bijection. In general situation, Lemma \ref{a24Feb13}.(1) shows that the map $\pil'$ is an injection. In general, there is no connection between maximal left denominator sets of  a ring and  its factor  ring.

Before giving the proof of Theorem \ref{24Feb13} we need some results.

\begin{lemma}\label{a24Feb13}
Let $R$ be a ring, $\gll = \gll_R$ be its left localization radical and $\pil : R\ra R/ \gll$, $r\mapsto \br = r+\gll$. Then
\begin{enumerate}
\item The map $\pil': \maxDen_l(R)\ra \maxDen_l(R/  \gll )$, $ S\mapsto \pil (S)$, is an injection such that $S^{-1}R\simeq \pil (S)^{-1} (R/\gll )$ and $\ass_{R/ \gll}(\pil (S))= \ass_R(S)/ \gll$.
\item The map $\pil'': \assmaxDen_l(R)\ra \assmaxDen_l(R/ \gll )$, $ \ga \mapsto \ga / \gll$, is an injection.
\item Let $T\in \maxDen_l(R/ \gll )$. Then $T\in \im (\pil')$ iff
 $\ass_{R/ \gll}(T)= \ass_R(S)/ \gll$ for some $S\in \maxDen_l(R)$.
 \item The map $\pil'$ if a bijection iff the map $\pil''$ is a bijection iff for every $\gb\in \maxDen_l(R/ \gll )$ there is $S\in \maxDen_l(R)$ such that for each element $b\in \pil^{-1}(\gb )$ there is an element $s\in S$ such that $ sb=0$.
\item If the map $\pil''$ is a bijection then its inverse $\pil'^{-1}$ is given by the rule $T\mapsto \pil'^{-1}(T)$.
\end{enumerate}
\end{lemma}

{\it Proof}. 1. Since $\gll \subseteq \ga = \ass (S)$ we have $\pil (S) \in \Den_l(R/\gll , \ga / \gll )$ and the map
\begin{equation}\label{pils}
\pilS : S^{-1}R\ra \pil (S)^{-1}R/ \gll, \;\; s^{-1}r\mapsto \pil^{-1}(s)\pil (r),
\end{equation}
is a ring isomorphism, by Lemma 3.2.(1), \cite{larglquot}. There is a commutative diagram of ring homomorphisms:
\begin{equation}\label{RRlS}
\xymatrix{ R\ar[r]^{\s_S}\ar[d]_{\pil}  & S^{-1}R\ar[d]^{\pilS} \\
 R/ \gll \ar[r]^{\overline{\s}_S}  & \pil (S)^{-1}R/ \gll  }
\end{equation}
where $\s_S(r) = \frac{r}{1}$ and $\overline{\s}_S(r+\gll ) = \frac{\pil (r)}{1}$. The ring $S^{-1}R$ is a left localization maximal ring. Let $G$ and $G'$ be the groups of units of the rings $S^{-1}R$ and  $\pil (S)^{-1}R/ \gll$ respectively. By Theorem \ref{15Nov10}.(3),
$$S= \s_S^{-1}(G)\;\; {\rm and } \;\; S':=\overline{\s}_S^{-1}(G')\in \maxDen_l(R/ \gll ).$$
By the commutativity of the diagram (\ref{RRlS}), we have
\begin{equation}\label{SpLS}
S= \pil^{-1} (S').
\end{equation}
Then, by the surjectivity of the map $\pil$, $\pil (S) = S'\in \maxDen_l(R/ \gll )$. So, the map $\pil'$ is a well-defined map which is an injection, by statement 2.

2. Every element of $\assmaxDen_l(R)$ contains $\gll$, hence $\pil'$ is an injection.

3. Statement 3 follows from statement 2 and Corollary \ref{d4Jan13}.(2).

4. By Corollary \ref{d4Jan13}.(2) and statements 1 and 2, the map $\pil'$ is a bijection/surjection iff the map $\pil''$ is a bijection/surjection. Let LS stands for the statement after the second `iff' in statement 4 and let $T\in \maxDen_l(R/ \gll)$ with $\gb = \ass (T) $. If $\pil'$ is a bijection then LS holds by statement 1 and the inclusion $\gll \subseteq \ass (S)$ for all $S\in \maxDen_l(R)$. If LS holds then clearly $\pil^{-1}(\gb )\subseteq \ass (S)$, and so $\gb \subseteq \ass (S) / \gll$. Then necessarily $\gb = \ass (S) / \gll$, by Proposition \ref{b27Nov12}, and so $T= \pil (S)$, by statement 1 and Corollary \ref{d4Jan13}.(2). This means that the map $\pil'$ is a surjection hence a bijection.

5. Statement 5 was proved in the proof of statement 1, see (\ref{SpLS}). $\Box $


\begin{lemma}\label{b24Feb13}
Let $R$ be a ring, $\gll = \gll_R$ be its left localization radical  and $\Min (R, \gll )$ be the set of  minimal primes over $\gll$. Suppose that $\maxDen_l(R)= \{ S_1, \ldots , S_n\}$ is a finite set and the ideals $\ga_i := \ass (S_i)$, $i=1, \ldots , n$, are prime. Then
\begin{enumerate}
\item The map $\maxDen_l(R)\ra \Min (R, \gll )$, $ S\mapsto \ass (S)$, is a bijection.
\item $\gll = \bigcap_{\gp \in\Min (R, \gll )}\gp$.
\end{enumerate}
\end{lemma}

{\it Proof}.  If $n=1$ then $\gll = \ga_1$ is a prime ideal and so statements 1 and 2  obviously hold.

So, let $n\geq 2$. The ideals $\ga_1, \ldots , \ga_n$ are distinct prime ideals with $\gll = \cap_{i=1}^n \ga_i$. Therefore, the prime ideals $\ga_1/ \gll, \ldots , \ga_n / \gll$ of the ring $R/ \gll$ are distinct with zero intersection. By Proposition \ref{A22Jan13}, the ring $R/ \gll$ is a semisimple ring with $\Min (R/ \gll )=\{ \ga_1/ \gll, \ldots , \ga_n / \gll )$. Hence, $\Min (R, \gll ) = \{ \ga_1, \ldots , \ga_n\}$ with $\gll = \bigcap_{i=1}^n \ga_n$, i.e. statement 2 holds. Now, statement 1 is obvious (by Corollary \ref{d4Jan13}.(2)). $\Box $

$\noindent $

{\bf Proof of Theorem \ref{24Feb13}}. $(1\Rightarrow 2)$  This implication follows from Theorem \ref{31Dec12} and Lemma \ref{a24Feb13}.(1).

$(2\Rightarrow 1)$  By (\ref{llRseq}), there is a ring monomorphism $R/ \gll \ra \prod_{S\in \maxDen_l(R)}S^{-1}R$. By Theorem \ref{22Jan13}, the ring $R$ has a semisimple left  quotient ring. Then, by Theorem \ref{22Jan13} and Lemma \ref{a24Feb13}.(1, 5), the statement (i) holds. By Goldie's Theorem  the ring $R$ is a semiprime left Goldie  ring.

$(1\Leftrightarrow 3)$ Notice that the map $\Min (R. \gll ) \ra \Min (R/ \gll )$, $\gp\mapsto \gp / \gll $, is a bijection. By Theorem \ref{1Jan13}, statement 1 is equivalent to the statements (a)'--(d)' (we have added `prime' to distinguish these conditions from the conditions (a)--(d) of statement 3):

(a)' $R/ \gll$ is a semiprime ring,

(b)' $|\Min (R/ \gll )|<\infty$,

(c)' For each $\gp / \gll \in \Min (R/ \gll )$, $S_{\gp / \gll }:=\{ c+\gll \in R/ \gll \, | \, c+\gp \in \CCRp \}\in \Den_l(R/ \gll , \gp / \gll )$, and

(d)' For each $\gp / \gll \in \Min (R/ \gll )$, the ring $S_{\gp / \gll }^{-1}R/ \gll $ is a simple left Artinian ring.

The conditions (a)' and (b)' are equivalent to the conditions (a) and (b) of statement 3. Now, the implication  $(1\Rightarrow 3)$  follows: the statement (i) and Lemma \ref{a24Feb13}.(4) imply the condition (e); then the conditions (c) and (d) follow from the conditions (c)' and (d)' by using Proposition \ref{a20Jan13}.

The implication  $(3\Rightarrow 1)$  follows from Theorem \ref{22Jan13}: the conditions (a)--(c) imply that there is a ring homomorphism
$$R/ \gll \ra \prod_{\gp\ \in \Min (R, \gll )}S_\gp^{-1}R.$$
By the condition (d), the rings  $S_\gp^{-1}R$ are simple left Artinian rings. Then, by Theorem \ref{22Jan13}, the ring $R/ \gll$ has a semisimple left quotient ring.

$(1\Rightarrow 4)$ The conditions (a) and (b) of statement 4 are precisely the conditions (a) and (b) of statement 3.  The condition (c) of statement 4 follows from the conditions (c) and (d) of statement 3 and Goldie's Theorem (since $Q(R/ \gp )\simeq S_\gp^{-1}R$ is a simple left Artinian ring, by statement 3).

$(4\Rightarrow 1)$ By the condition (c) and Goldie's Theorem, the left quotient rings $Q( R/ \gp )$ are  simple left Artinian rings. Then there is a  chain of ring homomorphisms
$$ R/ \gll = R/ \cap_{\gp \in \Min (R, \gp )}\gp\ra \prod_{\gp \in \Min (R, \gp )}R/ \gp \ra \prod_{\gp \in \Min (R, \gp )}Q(R/ \gp ).$$
By Theorem \ref{22Jan13}, the ring $R/ \gll$ has a semisimple left quotient ring. This finishes the proof of equivalence of statements 1--4. Now, statements (ii)--(iv) follow from the statement (i) and Lemma \ref{a24Feb13}.  $\Box $

$\noindent $

$\noindent $

$${\bf Acknowledgements}$$

 The work is partly supported by  the Royal Society  and EPSRC.

\small{

Department of Pure Mathematics

University of Sheffield

Hicks Building

Sheffield S3 7RH

UK

email: v.bavula@sheffield.ac.uk}

\end{document}